\documentclass{article}
\usepackage{amsmath}
\usepackage{amsthm}
\usepackage{amssymb}
\usepackage{euscript}

\usepackage[all]{xy}

\usepackage{diagrams}

\newenvironment{system}{%
\left \{ \begin{array}{@{\,}l@{\,}c@{\,}l@{\,}}%
  }%
  {%
  \end{array}\right\}
}
\newenvironment{myclaim}{\par\noindent\textrm{Claim}.~\it}{\par\noindent\relax}

\theoremstyle{plain}
\newtheorem{Theorem}{Theorem}[section]
\newtheorem{Proposition}[Theorem]{Proposition}
\newtheorem{Corollary}[Theorem]{Corollary}
\newtheorem{Lemma}[Theorem]{Lemma}

\theoremstyle{definition}
\newtheorem{Definition}[Theorem]{Definition}
\newtheorem{Example}[Theorem]{Example}
\newtheorem{Claim}[Theorem]{Claim}
\newtheorem{Fact}[Theorem]{Fact}

\theoremstyle{remark}

\renewcommand{\lefteqn}[2][0mm]{\makebox[#1][l]{$\displaystyle{#2}$}}

\newcommand{\timplies}{\;\textrm{implies}\;}
\newcommand{\tor}{\;\textrm{or}\;}
\newcommand{\tst}{\;\textrm{s.t.}\;}

\newcommand{\MN}[1]{\overline{#1}}
\newcommand{\ideal}[1]{i #1}
\newcommand{\filter}[1]{\uparrow\!{#1}}
\newcommand{\righta}[1][]{\mathrm{r}_{#1}}
\newcommand{\Cvrs}[2]{\EuScript{C}(#1;#2)}

\newcommand{\arrow}[1][]{\stackrel{#1}{\rightarrow}}
\newcommand{\spcon}[1][]{\mbox{$\displaystyle\bigwedge\!\!\!\!\!\!\!\bigwedge$}_{#1}}
\newcommand{\nec}[1][\hspace{0.5em}]{[#1]}
\newcommand{\pos}[1][\hspace{0.5em}]{\langle#1\rangle}

\newcommand{\Opens}{\mathcal{O}_{\!f}}

\newcommand{\Alg}[1]{\mathcal{#1}}
\newcommand{\T}{\mathbb{T}}

\newcommand{\id}{\mathtt{id}}
\newcommand{\prj}[1][]{{\mathtt{pr}}^{#1}} 

\newcommand{\set}[1]{\{\,#1\,\}}
\newcommand{\quotient}{\!/\!\!\sim}
\newcommand{\ppowerset}{\mathcal{P}_{\!\!+}}
\newcommand{\reach}[1]{\overline{#1}}
\newcommand{\myjoin}[1][]{\textstyle{\bigvee_{\!#1}}}
\newcommand{\nnumbers}{\mathbb{N}}

\newcommand{\wi}{l}
\newcommand{\wh}{m}
\newcommand{\nh}{h}
\newcommand{\nii}{i}

\newcommand{\mydiagram}[2][6em]{
  \mbox{$\begin{array}{@{\hspace{0mm}}c@{\hspace{0mm}}}
      \xy\xygraph{!{<#1,0cm>:<0cm,#1>::}#2}\endxy
    \end{array}$}
  }

\newgraphescape{s}[6]{
  []*+{#1}="1"
  (
  [d(#5)]*+{#2}="2"
  [r(#6)]*+{#4}="4",
  [r(#6)]*+{#3}="3"
  )
  }

\newgraphescape{a}[4]{
  "1":"2"#1
  "1":"3"#2
  "2":"4"#3
  "3":"4"#4
}

\author{Luigi Santocanale\\[2mm]%
  LIF-CMI Marseille\\[2mm]%
  \texttt{luigi.santocanale@cmi.univ-mrs.fr}%
}
\title{Completions of $\mu$-Algebras}
\date{\today}
\begin{document}

\maketitle

\begin{abstract}

A $\mu$-algebra is a model of a first order theory that is an
extension of the theory of bounded lattices, that comes with pairs of
terms $(f,\mu_{x}.f)$ where $\mu_{x}.f$ is axiomatized as the least
prefixed point of $f$, whose axioms are equations or equational
implications.

Standard $\mu$-algebras are complete meaning that their lattice reduct
is a complete lattice.  We prove that any non trivial quasivariety of
$\mu$-algebras contains a $\mu$-algebra that has no embedding into a
complete $\mu$-algebra.

We focus then on modal $\mu$-algebras, i.e. algebraic models of the
propositional modal $\mu$-calculus. We prove that free modal
$\mu$-algebras satisfy a condition -- reminiscent of Whitman's
condition for free lattices -- which allows us to prove that (i) modal
operators are adjoints on free modal $\mu$-algebras, (ii) least
prefixed points of $\Sigma_{1}$-operations satisfy the constructive
relation $\mu_{x}.f = \bigvee_{n \geq 0} f^{n}(\bot)$.  These
properties imply the following statement: {\em the MacNeille-Dedekind
  completion of a free modal $\mu$-algebra is a complete modal
  $\mu$-algebra and moreover the canonical embedding preserves all the
  operations in the class $Comp(\Sigma_{1},\Pi_{1})$ of the
  fixed point alternation hierarchy.}

\end{abstract}


\section*{Introduction}

When $L$ is a complete lattice, the least fixed point $\mu_{x}.f$ of a
monotone function $f: L \rTo L$ enjoys a remarkable property. We like
to say that the least fixed point is \emph{constructive}: the
equality
\begin{align}
  \label{eq:approximants}
  \mu_{x}.f
  & = \bigvee_{\alpha \in Ord} f^{\alpha}(\bot)
\end{align}
holds and provides a method to construct $\mu_{x}.f$ from the bottom
of the lattice.  The expressions $f^{\alpha}(\bot)$, indexed by
ordinals, are commonly called the \emph{approximants} of $\mu_{x}.f$.
They are defined by transfinite induction as expected:
$f^{0}(\bot)=\bot$, $f^{\alpha + 1} = f(f^{\alpha}(\bot))$, and
$f^{\alpha}(\bot) = \bigvee_{\beta < \alpha} f^{\beta}(\bot)$ for a
limit ordinal $\alpha$.
                                
A careful reading of Tarski's original fixpoint theorem
\cite{tarski55} reveals that the completeness assumption is not needed
for $f$ to have a fixed point.  If $L$ is merely a poset a \emph{least
  prefixed point} of a monotone $f$ is an element $\mu_{x}.f\in L$
satisfying
\begin{align}
  \label{eq:fix}
  f(\,\mu_{x}.f\,) & \leq \mu_{x}.f \,,\\
  \label{eq:park}
  f(y) \leq y \;\;\;\;\;\; \Rightarrow 
  \;\;\;\;\;\;\mu_{x}.f & \leq y\,.
\end{align}
Tarski's theorem can be rephrased by saying that the least prefixed
point, whenever it exists, is also a fixed point, hence it is the least
fixed point. The two notions are similar and coincide on complete
lattices.  Properties \eqref{eq:fix} and \eqref{eq:park} -- the latter
known as the Park induction rule \cite{park,esik} -- provide a natural
axiomatization by equations and equational implications of least
fixed points, provided that an order relation definable by equations
is given.  They have been used often to axiomatize concrete
mathematical objects where implicit or explicit fixed points are at
work: relational algebras with transitive closure \cite{NT77}, regular
languages \cite{kozen:kleenealgebras}, powersets of Kripke frames
\cite{segerberg,kozen}.  

Many considerations induce to study 
classes of models of axioms \eqref{eq:fix} and \eqref{eq:park}. For
example, model theory suggests that models of theories axiomatized by
equational implications are preferable to models that are complete
lattices: the former build up a quasivariety, colimits exist, free
models exist, etc.  The goal of this paper is to compare the models of
theories where the least fixed points are defined by means of
\eqref{eq:fix} and \eqref{eq:park} -- we shall refer to them as
$\mu$-algebras -- with a more restricted class of models, the
standard, concrete, or complete models. These are the models whose
underlying lattice is complete and where least fixed points are
constructive.  Mathematically, the comparison amounts to asking
whether a model can be embedded into a complete one.

Despite the difference in the respective lengths, the paper is divided
into two parts. In the first part we show that almost never
$\mu$-algebras are completable. That is, within a non trivial fixed
quasi-variety of $\mu$-algebras, we construct a $\mu$-algebra that has
no embedding into a complete one.  The second part of the paper is
devoted to studying free modal $\mu$-algebras. Modal $\mu$-algebras
are algebraic models of the modal $\mu$-calculus \cite{kozen}; by the
completeness theorem w.r.t. the class of Kripke frames
\cite{kozen,wal} we already know that free modal $\mu$-algebras are
completable.  We pursue an algebraic understanding of this fact, which
eventually will provide us with some algebraic interpretation of the
completeness theorem.  Our analysis of free $\mu$-algebras, which
never takes the completeness theorem as granted, can be synthesized as
follows. We observe first a phenomenon that we classify as
``definability of adjoints''.  Using adjoints and their
generalizations, $\Opens$-adjoints, we argue that a restricted class
of least fixed points are constructive on free modal $\mu$-algebras.
This means that many relations like \eqref{eq:approximants} hold on a
free modal $\mu$-algebra even if this is presumably not complete. A
detour through least solutions of systems of equations allow us to
extend the class of constructive operations on free modal
$\mu$-algebras to include all the $\Sigma_{1}$-operations.  It is
easily argued that constructiveness is an essential property for an
embedding into a complete $\mu$-algebra to exist. Indeed, the outcome
of our analysis is the following result: the MacNeille-Dedekind
completion of a free modal $\mu$-algebra is a complete modal
$\mu$-algebra and the canonical embedding preserves all the operations
in the class $\Sigma_{1}$ of the fixed point alternation hierarchy.
The result is easily extended to the class $Comp(\Sigma_{1},\Pi_{1})$
of the alternation hierarchy.

\tableofcontents


\section{Notation, background}

\subsection*{Cartesian structure}

If $X$ is a finite set, then $L^{X}$ will denote the power of $L$ by
$X$, i.e. $L^{X} = \prod_{x \in X} L$.  Projection functions from a
product to one of its factors will be denoted by $\prj$, with the
necessary sub/superscripts. For example if $x \in X$, then
$\prj[X]_{x} : L^{X} \rTo L$ denotes the projection taking a vector $v
\in L^{X}$ to $v(x)$.

If $f_{x}:L \rTo M$, $x \in X$, is a collection of functions, then
shall use the notation $\langle f_{x} \rangle_{x \in X}: L \rTo M^{X}$
for the unique function $f : L \rTo M^{X}$ such that $f_{x} = \prj_{x}
\circ f$.

\subsection*{Parametrized fixed points}
\label{subsec:parameters}
Let $f: P \times Q \rTo P$ be a monotone function.  For $q \in Q$, we
use the notation $f_{q} : P \rTo P$ for the monotone function sending
$p$ to $f(p,q)$.
\begin{Lemma}
  Suppose that for each $q
  \in Q$ the least prefixed point $\mu_{x}.f_{q}$ of $f_{q} : P \rTo P$
  exists. Then the correspondence $\mu_{x}.f : Q \rTo P$ sending $q$ to
  $\mu_{x}.f_{q}$ is monotone.
\end{Lemma}
The notation $f: P^{x}\times Q^{y} \rTo M$ will be used here to mean
that $f$ is considered as a function of two variables $x \in P$ and $y
\in Q$. For example, if $f: P^{x}\times P^{y} \rTo P$, then we shall
write $\mu_{y}.f$ or $\mu_{y}.f(x,y)$ for unary function sending $p
\in P$ to the least prefixed point of the unary function $y \rMapsto
f(p,y)$.

\subsection*{The Beki\v{c} property}

We shall often make use of the Beki\v{c} property which usually is
stated as an identity between existing (least pre-)fixed points.  We
shall also be concerned with existence of least prefixed points, hence
we need a stronger form of the property which emphasizes this issue as
well.
\begin{Proposition}
  Let $P,Q$ be posets, let $f : P^{x} \times Q^{y} \rTo P$ and $g : P
  \times Q \rTo \ Q$ be monotone functions, and suppose that for each
  $q \in Q$ the least prefixed point $\mu_{x}.f_{q}$ of $f_{q} : P
  \rTo P$ exists.  Consider the monotone functions
  \begin{align*}
    \langle f,g\rangle & : P \times Q \rTo P\times Q\,,\\
    \langle \mu_{x}.f \circ \prj_{Q}, g \rangle & : P\times Q
    \rTo P\times Q\,.
  \end{align*}
  A least prefixed point $(\mu_{1},\mu_{2})$ of $\langle f,g\rangle$
  exists if and only if the least prefixed point $(\mu_{3},\mu_{4})$ of
  $\langle \mu_{x}.f \circ \prj_{Q}, g \rangle$ exists, and if any of
  them exists then $(\mu_{1},\mu_{2}) = (\mu_{3},\mu_{4})$.
\end{Proposition}

\begin{Proposition}
  Let $P,Q$ be posets and  $f : Q \rTo P$, $g : P \times
  Q \rTo \ Q$ be monotone functions.
  Consider the monotone functions
  \begin{align*}
    \langle f \circ \prj_{Q}, g \rangle
    & :  P \times Q \rTo  P \times Q \,, \\
    g \circ \langle f,\id_{Q} \rangle
    & : Q \rTo P \times Q \rTo Q
    \,.
  \end{align*}
  A least prefixed point $(\mu_{1},\mu_{2}) \in P \times Q$ of
  $\langle f \circ \prj_{Q}, g \rangle$ exists if and only if a least
  prefixed point $\mu_{3} \in Q$ of $g \circ \langle f,\id_{Q}
  \rangle$ exists, and they are determined by each other as folllows:
  \begin{align*}
    (\mu_{1},\mu_{2}) & = (f(\mu_{3}),\mu_{3})\,,
     & \mu_{3} & = \mu_{2}  \,.
  \end{align*}
\end{Proposition}
A proof of the propositions appears in the extended version of
\cite{san:paritygames} within the more general context of initial
algebras of functors.


\section{$\mu$-algebras are not completable}

\subsection*{$\mu$-theories and $\mu$-algebras}

Our first goal is to setup a generic logical framework
within which to develop a theory of ordered algebras with least
fixed point operators.  Analogous frameworks \cite{bloomesik,AN} can
be coded within this framework.

\begin{Definition}
  A \emph{$\mu$-theory} is a first order theory with the following
  properties:
  \begin{itemize}
  \item it is an extension of the theory of bounded lattices,
  \item it comes with fixed point pairs, that is, pairs of terms
    $(f,\mu_{x}.f)$ axiomatized by \eqref{eq:fix} and \eqref{eq:park}
    so that (the interpretation of) $f$ is an order preserving
    operation in the variable $x$, and (the interpretation of)
    $\mu_{x}.f$ is a least prefixed point of $f$,
  \item its axioms are either equations or equational implications.
  \end{itemize}
  A \emph{$\mu$-algebra} is model of a fixed $\mu$-theory.  A
  $\mu$-algebra is \emph{complete} if its lattice reduct is a complete
  lattice.
\end{Definition}
The notion of a morphism of $\mu$-algebras is standard from model
theory: a function $g : A \rTo B$ between the underlying sets of
$\mu$-algebras $\Alg{A}$ and $\Alg{B}$ is a morphism if it preserves
the interpretation of all the terms of the $\mu$-theory.

Let $f$ be a term of a $\mu$-theory and let $X$ be its set of free
variables. For a $\mu$-algebra $\Alg{A}$, we shall overload the
notation and write $f: \Alg{A}^{X} \rTo \Alg{A}$ for the
interpretation of $f$ on $\Alg{A}$. If $f$ is part of a fixed point
pair $(f,\mu_{x}.f)$, so that $X$ is the disjoint union of $\set{x}$
and $Y$, and $v \in \Alg{A}^{Y}$, then we use the notation $f_{v} :
\Alg{A} \rTo \Alg{A}$ consistenlty with what exposed in Section
\ref{subsec:parameters} and say that $f_{v}: \Alg{A} \rTo \Alg{A}$ is
a fixed point polynomial.  To simplify the notation, we shall also
omit the subscript $v$ and say that $f: \Alg{A} \rTo \Alg{A}$ is a
fixed point polynomial.

\subsection*{A $\mu$-algebra with no complete extension}

Our goal is to assess relations between $\mu$-algebras and complete
$\mu$-algebras and to understand when a $\mu$-algebra
embeds into a complete one.  Contrarily to what happens for several
algebraic structures related to logic (Boolean algebras, modal
algebras K, Heyting algebras, quantales), we show next that this is
not always possible for $\mu$-algebras.

\begin{Example}
  \label{ex:nonembeddable}
  Choose a $\mu$-algebra $\Alg{A}$ and a fixed point polynomial $f :
  \Alg{A} \rTo \Alg{A}$ for which the chain of finite approximants
  $$
  \bot < f(\bot) < f^{2}(\bot) < \ldots < f^{n}(\bot) < \ldots 
  $$
  is infinite.  Define the infinite sequences $\phi_{n}$ by
  \begin{align*}
    \phi_{n} & = (\;\underbrace{\bot,\ldots,\bot}_{n-\textrm{times}},
    \bot,f(\bot),f^{2}(\bot),\ldots \;)\,,\;\;n \geq 0\,,
  \end{align*}
  and consider them as elements of the product algebra
  $\Alg{A}^{\omega}$.  Since $f$ is computed pointwise, observe
  that $f(\phi_{n})$ is equal to $\phi_{n -1}$ for all but a finite
  number of coordinates.
  
  Define the equivalence relation $\sim$ on $\Alg{A}^{\omega}$ by
  saying that two infinite sequences are equivalent if they coincide
  in all but a finite number of coordinates. The quotient
  $\Alg{A}^{\omega}\quotient$ is a  reduced product of $\Alg{A}$
  and all the equations and equational implications that hold in
  $\Alg{A}$ hold in $\Alg{A}^{\omega}\quotient$ as well, cf.
  \cite[chapter 6]{changkeisler}. In particular,
  $\Alg{A}^{\omega}\quotient$ is  a $\mu$-algebra in the same
  quasivariety as $\Alg{A}$.
  
  Denote by $\bar{\phi}_{n}$ the equivalence class of $\phi_{n}$ and
  recall that in $\Alg{A}^{\omega}\quotient$ the least fixed point
  $\mu_{x}.f$ is simply the equivalence class of the infinite sequence
  with constant value $\mu_{x}.f$. The relations
  \begin{align}
    \label{eq:wrongconf}%
    f(\bar{\phi}_{n}) & \leq
    \bar{\phi}_{n-1},\,\;\;  n \geq 1 &
    \mu_{x}.f &\not\leq \bar{\phi}_{0}%
  \end{align}
  hold in $\Alg{A}^{\omega}\quotient$ and we claim that a
  configuration such as the one described by \eqref{eq:wrongconf} is
  not compatible with $\Alg{A}^{\omega}\quotient$ being complete. If
  $\bigwedge_{n \geq 0} \bar{\phi}_{n}$ exists then
  $$
  f(\bigwedge_{n \geq 0} \bar{\phi}_{n}) \;\leq\; f(\bar{\phi}_{n +
    1}) \; \leq \bar{\phi}_{n}\,,
  $$
  for all $n \geq 0$, and therefore $f(\bigwedge_{n \geq 0}
  \bar{\phi}_{n})\leq \bigwedge_{n \geq 0} \bar{\phi}_{n}$. Then
  $\mu_{x}.f \leq \bigwedge_{n \geq 0} \bar{\phi}_{n} \leq
  \bar{\phi}_{0}$ gives a contradiction.
  
  Finally observe that such a configuration is preserved by any
  extension of $\Alg{A}^{\omega}\quotient$, and therefore this
  $\mu$-algebra has no complete extension. \qed
\end{Example}

It can be observed that in the $\mu$-algebra
$\Alg{A}^{\omega}\quotient$ the stronger relations $f(\bar{\phi}_{n})
= \bar{\phi}_{n-1}$ hold and moreover $\phi_{n + 1} \leq \phi_{n}$.
These stronger relations, which are not needed to prove that
$\Alg{A}^{\omega}\quotient$ has no completion,
conceal the original idea behind the impossibility proof. This amounts
to the
construction of a chain of approximants indexed by natural numbers
with the reverse order.

We say that a $\mu$-theory (or the quasivariety of the $\mu$-algebras)
is \emph{non-trivial} if we can find a $\mu$-algebra $\Alg{A}$ and a
fixed point polynomial $f$ for which its finite approximants are all
distinct.  If this is not possible, then for each fixed point pair
$(f,\mu_{x}.f)$ some equation of the form $\mu_{x}.f = f^{n}(\bot)$
holds, showing that all the least fixed point are superflous. We
collect these observations in a Theorem.
\begin{Theorem}
  Any non-trivial quasivariety of $\mu$-algebras contains a
  $\mu$-algebra which does not admit an embedding into a complete
  $\mu$-algebra.
\end{Theorem}
For simple $\mu$-theories, if a
$\mu$-algebra has no configuration such as \eqref{eq:wrongconf}, then
the principal filter embedding is a morphism of $\mu$-algebras.

\begin{Example}
  Consider a $\mu$-theory $\T$ with just a unary function symbol $f$
  and a constant $\mu_{x}.f$ in addition to the signature of bounded
  lattices.  The axioms of $\T$ are those of bounded lattices,
  additional equations, and the
  fixed point axioms \eqref{eq:fix} and \eqref{eq:park} for the unique
  fixed point pair $(f(x),\mu_{x}.f)$.
  
  Let $\Alg{A}$ be a $\mu$-algebra with no configuration such as
  \eqref{eq:wrongconf} and let $\mathcal{F}(\mathcal{A})$ be the
  standard algebra of filters of $\Alg{A}$ in the signature of $\T$.
  Then all the equations of $\T$ holds in $\mathcal{F}(\mathcal{A})$
  and all we need to observe is that the principal filter $\filter
  \mu_{x}.f$ is the least fixed point of the extension of $f$ to
  $\mathcal{F}(\mathcal{A})$. Considering that the order in
  $\mathcal{F}(\mathcal{A})$ is reverse inclusion, we need to verify
  that $\mu_{x}.f$ is below any element of an arbitrary filter $F$
  such that $F \subseteq f(F)$.  Recalling that
  \begin{align*}
    f(F) & = \set{  x\,|\,\exists y \in F f(y) \leq x\,}\,,
  \end{align*}
  if $\phi_{0} \in F$ then we can construct a sequence
  $\set{\phi_{n}}_{n \geq 0}$ such that $f(\phi_{n +1}) \leq \phi_{n}$
  for $n \geq 0$. Since $\mathcal{A}$ lacks a configuration such as
  \eqref{eq:wrongconf}, we deduce $\mu_{x}.f \leq \phi_{0}$.  \qed
\end{Example}
Continuity of the algebra of ideals over a lattice is a major obstacle
to exploit such a construction for completions of $\mu$-algebras. For
example, the principal filter embedding becomes useless for
$\mu$-theories where greatest fixed points are also an issue.  Other
conditions are needed to ensure that a $\mu$-algebra has an embedding
into a complete $\mu$-algebra.


\section[Completions for free modal $\mu$-algebras, overview]{Completions for free modal $\mu$-algebras,\\\mbox{\hspace{85mm}} overview}

Recall that a free $\mu$-algebra embeds into a complete one if and
only if the class of complete $\mu$-algebras generates the class of
all $\mu$-algebras. If we adopt the perspective of algebraic logic,
the statement that free $\mu$-algebras embed into complete ones
amounts to a completeness theorem for the logic with respect to the
semantics of all complete models.

It is often the case that free $\mu$-algebras embed into complete
ones, for example free $\mu$-lattices \cite{mulat} and free modal
$\mu$-algebras, i.e.  Lindenbaum algebras for the propositional modal
$\mu$-calculus \cite{kozen}.  The rest of this paper will be concerned
with studying \emph{free modal $\mu$-algebras}. We present here their
$\mu$-theory, i.e. the theory of modal $\mu$-algebras. The terms of
the theory are generated according to the  grammar:
\begin{align*}
  t & = p\,| 
  \,x\,|\,\top \,| \,t_{1} \land t_{2} \,
  |\,\neg t \,| \,\pos[\sigma] t\,|\,\mu_{x}.t\,,
\end{align*}
where $\sigma$ ranges on a finite set of actions $Act$ and the
fixed point generation rule applies only when the variable $x$
occurs under an even number of negations.  The reader has surely
recognized the framework of multimodal algebras, in addition to which
we have least fixed points.  Accordingly, the axioms of the theory are
those of multimodal algebras K as well as \eqref{eq:fix} and
\eqref{eq:park} for the fixed point pairs $(t,\mu_{x}.t)$.  In the
grammar we have distinguished a generator $p$ from a variable $x$.
This will be useful when considering the interpretation of terms as
operations on free modal $\mu$-algebras, where the generators become
operations.  This kind of term generation is standard from fixed point
theory \cite{niwinski2}, but it is also possible to code these terms
as terms generated from an infinite signature using substitution
only \cite{nelson}. Finally, it can be shown that modal $\mu$-algebras
form a variety of algebras \cite{san:eqfixedpoints}.

The completeness results for the propositional modal $\mu$-calculus
\cite{kozen,wal} paired with the small Kripke model property
\cite{emersonstrett} imply that \emph{a free modal $\mu$-algebra has
  an embedding into} an infinite product of finite modal
$\mu$-algebras. This infinite product is of course \emph{a complete
  lattice}.  In the rest of the paper we shall prove a weaker
embedding result concerning $\Sigma_{1}$-terms and
$\Sigma_{1}$-operations. $\Sigma_{1}$-terms are defined by the
grammar:
\begin{align}
 \label{eq:sigma1}
  t & \,= \;x \,|\, p\,|\,\neg p \,| \,\top \,| \,t \land t \,| \, \bot
  \,|\, t \vee t \,|\, \pos[\sigma] t\,|\,\nec[\sigma]t \,|\,
  \mu_{x}.t\,,
\end{align}
and $f : \Alg{A}^{X} \rTo \Alg{A}$ is a $\Sigma_{1}$-operation if it
is the interpretation of a $\Sigma_{1}$-term.  Observe that the
fixed point formation rule is no
longer constrained in the above grammar.  
By duality, the greatest fixed point $\nu_{x}.f(x,y)$ of an operation
$f(x,y)$ is definable in the given signature: $\nu_{x}.f(x,y) = \neg
\mu_{x}.\neg f( \neg x ,y)$.  The class of $\Pi_{1}$-terms is then
defined as above with the exception that least fixed point formation
is replaced by greatest fixed point formation.  The class of
$Comp(\Sigma_{1},\Pi_{1})$-operations is obtained by composing in all
the possible ways operations in the classes $\Sigma_{1}$ and
$\Pi_{1}$.  The reader is invited to consult \cite[chapter 8]{AN} for
an exposition of the full fixed point alternation hierarchy. Our
result can be stated as follows:
\begin{Theorem}
  \label{theo:maintarget}
  Let $\mathcal{F}$ be a free modal $\mu$-algebra. There exists a
  complete modal algebra $\MN{\mathcal{F}}$ and an injective
  morphism of Boolean modal algebras $\ideal : \mathcal{F} \rTo
  \MN{\mathcal{F}}$ which preserves all the
  $Comp(\Sigma_{1},\Pi_{1})$-operations of the algebra $\mathcal{F}$.
\end{Theorem}
With respect to \cite{wal}, where algorithmic and game-theoretic ideas
as well as tableaux manipulations are the main tools, we shall use 
purely algebraic and order theoretic tools. Under some respect, our
work can also be understood as an effort to translate ideas from
\cite{kozen,wal} into an algebraic and order theoretic framework.

We sketch in the rest of the section the strategy followed to prove
Theorem \ref{theo:maintarget}.  The algebra $\MN{\mathcal{F}}$ is the
MacNeille-Dedekind completion of $\mathcal{F}$.  For our goals, we
recall that if $L$ is a Boolean algebra, then $\MN{L}$ is a Boolean
algebra as well, see  \cite[Chapter V, Theorem 27]{birkhoff}.
Recall that an order preserving $f : L \rTo M$ is a left adjoint if
there exists $g : M \rTo L$ (the right adjoint) such that $f(x) \leq
y$ if and only if $x \leq g(y)$, for all $x \in L$ and $y \in M$.  For
our goals, we also need the following statement:
\begin{Lemma}
  \label{lemma:adjmcneille}
  Let $L$ be a lattice and $\MN{L}$ be its MacNeille-Dedekind
  completion.  A left adjoint $f : L \rTo L$ has an extension --
  necessarily unique -- to a left adjoint $f^{\vee} : \MN{L} \rTo
  \MN{L}$.
\end{Lemma}
Using the notation of \cite{harding}, if $g$ is right adjoint to $f$,
then $g^{\wedge}$ is right adjoint to $f^{\vee}$.  A first step
towards our main result will be to prove:
\begin{Claim}
  \label{claim:uno}
  The  modal operators $\pos[\sigma]$ of a free modal
  $\mu$-algebra are left adjoints.
\end{Claim}
Using Lemma \ref{lemma:adjmcneille} and Claim \ref{claim:uno} we can
state:
\begin{Proposition}
  The MacNeille-Dedekind completion $\MN{\mathcal{F}}$ of a free modal
  $\mu$-algebra is a multi-modal algebra K and the
  principal ideal embedding is a morphism of multi-modal algebras.%
  \footnote{The same statement holds if we replace ``free modal
    $\mu$-algebra'' with ``free multi-modal algebra K''.}
\end{Proposition}
Since $\MN{\mathcal{F}}$ is a complete lattice, it is a complete modal
$\mu$-algebra, and therefore we are also interested in preservation of
fixed points. To this goal we shall use the following Lemma:
\begin{Lemma}
  \label{lemma:pres}
  Let $\Alg{A}$ be a $\mu$-algebra, $i : \Alg{A} \rTo \MN{\Alg{A}}$
  its MacNeille-Dedekind completion, and $f_{v}$ a fixed point polynomial.
  Suppose that
  \begin{itemize}
  \item $f_{v}$ is preserved by $i$, that is, $i(f_{v}(x)) =
    f_{i(v)}(i(x))$,
  \item $\mu_{x}.f_{v}$ is constructive: $\mu_{x}.f_{v} =
    \bigvee_{\alpha \in Ord} f_{v}^{\alpha}(\bot)$.  
  \end{itemize}
  Then the least fixed point $\mu_{x}.f_{v}$ is preserved:
  $i(\mu_{x}.f_{v}) = \mu_{x}.f_{i(v)}$.
\end{Lemma}
\begin{proof}
  Observe first that
  $$
  f_{i(v)}(i( \mu_{x}.f_{v}))\; = \;i(f_{v}( \mu_{x}.f_{v}) ) \;=
  \;i(\mu_{x}.f_{v})\,,$$
  from which we deduce $\mu_{x}.f_{i(v)} \leq
  i(\mu_{x}.f_{v})$.  For the converse we argue that approximants are
  preserved using continuity of the embedding of a lattice into its
  MacNeille-Dedekind completion.  We have that $f_{v}^{0}(\bot)$ is
  preserved since $i$ preserves the bottom, and
  $f_{v}^{\alpha+1}(\bot)$ is preserved since $i$ preserves $f$. For a
  limit ordinal $\alpha$, suppose that $i$ preserves $f^{\beta}(\bot)$
  for $\beta < \alpha$. Then:
  \begin{align*}
    i(\bigvee_{\beta < \alpha} f_{v}^{\beta}(\bot)) 
    & =  \bigvee_{\beta < \alpha} i (f_{v}^{\beta}(\bot)) 
    =  \bigvee_{\beta < \alpha}  f_{i(v)}^{\beta}(\bot) \,,
  \end{align*} 
  since $i$ preserves all existing joins.  Consequently
  $i(\mu_{x}.f_{v}) = \bigvee_{\alpha \in Ord} f_{i(v)}^{\alpha}(\bot)$
  which clearly is below $\mu_{x}.f_{i(v)}$.  
\end{proof}
We shall prove that all the $\Sigma_{1}$-operations are preserved by
showing that all these functions are constructive:
\begin{Claim}
  \label{claim:due}
  Every fixed point $\Sigma_{1}$-polynomial
  $f_{v} : \mathcal{F} \rTo \mathcal{F}$ over a free modal
  $\mu$-algebra satisfies the constructive relation
  \begin{align}
    \label{eq:constr}
    \mu_{x}.f_{v} & = \bigvee_{n \geq 0} f_{v}^{n}(\bot)\,.
  \end{align}
\end{Claim}
Lemma \ref{lemma:pres} and Claim \ref{claim:due} imply that each
$\Sigma_{1}$-operation on a free modal $\mu$-algebra is preserved. 

A proper dualisation of the notions and results exposed so far can be
used to prove that $\Pi_{1}$-operations are preserved by the embedding
of a free modal $\mu$-algebra into its MacNeille-Dedekind completion.
Consequently, all the operations in the class
$Comp(\Sigma_{1},\Pi_{1})$ are preserved as well.

Finally, it should be observed that -- using the completeness of the
propositional modal $\mu$-calculus and the small model theorem
\cite{kozen,wal,emersonstrett} -- it is possible to directly argue
that every fixed point polynomial on a free modal $\mu$-algebra
satisfies the relation \eqref{eq:constr}. Hence the embedding of a
free modal $\mu$-algebra into its MacNeille-Dedekind completion is
indeed a morphism of modal $\mu$-algebras. On the other hand, it is
implicit from \cite{wal2} that the constructive relation
\eqref{eq:constr} is at the core of the completeness problem for the
modal $\mu$-calculus. There a proof-system is described which can
easily be proved complete once it is known that the relations
\eqref{eq:constr} hold in a free modal $\mu$-algebras.

We introduce the main property of free modal $\mu$-algebras in the
next section. Using this property we shall immediately be able to
prove Claim \ref{claim:uno}.  We introduce next the notion of a
$\Opens$-adjoint of finite type, using which we shall be able to prove
constructiveness of several monotone endofunctions on products of free
modal $\mu$-algebras, i.e. systems of equations on free modal
$\mu$-algebras. We shall prove Claim \ref{claim:due} at the end of the
paper, after devising a method to transport constructiveness from
systems of equations to operations.


\section{A property of free modal $\mu$-algebras}

In this section we prove that free modal $\mu$-algebras enjoy a
property similar to Whitman's condition for free lattices, cf.
\cite{whitman,freese}. In proof theory this sort of property is often
called a last rule and implies a cut-elimination theorem. We do not
know yet if this property, stated in the next Theorem, characterizes
free modal $\mu$-algebras.\footnote{It is easily argued that this
  property characterizes free modal algebras K among the finitely
  generated ones.}  However this property is  quite powerful and will
eventually allow us to prove Claim \ref{claim:uno} and Claim
\ref{claim:due}. 

We briefly recall the universal property of a modal $\mu$-algebra
$\mathcal{F}_{P}$ freely generated by a set $P$. Such a $\mu$-algebra
comes with a function $j : P \rTo \mathcal{F}_{P}$ such that for each
pair $(f,\Alg{A})$ -- where $\Alg{A}$ is a modal $\mu$-algebra and $f
: P \rTo \Alg{A}$ -- there exists a unique $\mu$-algebra morphism
$\tilde{f}: \mathcal{F}_{P} \rTo \Alg{A}$ such that $f = \tilde{f}
\circ j$. A generator in $\mathcal{F}_{P}$ is of the form $j(p)$ for
some $p \in P$.  It is easily argued that $j$ is injective (see the
end of this section) and therefore we shall abuse notation and
identify $P$ with its image $j(P)$.
\begin{Theorem}
  \label{theo:mainsec}
  \label{theo:mainproperty}
  Let $\mathcal{F}$ be a free modal $\mu$-algebra and $\Lambda$ be a
  finite set of literals (generators or negated generators).  The
  following implication holds in $\mathcal{F}$: if
 \begin{align*}
     \bigwedge \Lambda \land \bigwedge_{\sigma \in Act}
     (\,\nec[\sigma]x_{\sigma} \land \bigwedge_{y \in Y_{\sigma}}
     \pos[\sigma] y\,) & \leq \bot\,,
   \end{align*}
   then either $p,\neg p \in \Lambda$ for some generator $p$, or $
   x_{\sigma} \land y \leq \bot$ for some $\sigma \in Act$ and
   $y \in Y_{\sigma}$.
\end{Theorem}
We prove  first that:
\begin{Proposition}
  Let $\Alg{F}$ be a free modal algebra. The implication
  \begin{equation}
    \label{eq:mainresult}
    \begin{split}
      \bigwedge \Lambda 
      \;\land\;
      \bigwedge_{\sigma \in \Sigma} &
      (\;
      \nec[\sigma]\bigvee Y_{\sigma}%
      \;\land\;%
      \bigwedge_{y \in
        Y_{\sigma}} \pos[\sigma] y
      \;)%
      \leq \bot \\
      &\timplies \\
      \bigwedge \Lambda &\leq \bot \tor \exists \sigma \in
      \Sigma, y \in Y_{\sigma} \tst y \leq \bot
    \end{split}
  \end{equation}
  holds in $\Alg{F}$, where $\Lambda$ is a finite set of literals,
  $\Sigma \subseteq Act$, and, for each $\sigma \in \Sigma$,
  $Y_{\sigma}$ is a finite possibly empty set of elements of
  $\Alg{F}$.
\end{Proposition}
\begin{proof}
  Let $\Alg{A}$ be any modal algebra and suppose that for each $\sigma
  \in \Sigma$ we are given a set $Y_{\sigma}$ such that $y \not\leq
  \bot$ for each $y \in Y_{\sigma}$.  For each $\sigma \in \Sigma$ and
  $y \in Y_{\sigma}$ let $\chi_{\sigma}^{y}: \Alg{A} \rTo 2$ be
  morphism of Boolean algebras such that $\chi_{\sigma}^{y}(y) =
  \top$ (such a morphism exists by the prime filter theorem). Define
  \begin{align*}
    \chi_{\sigma}(z) 
    & = 
    \begin{cases}
      \bigvee_{y \in Y_{\sigma}} \chi_{\sigma}^{y}(z)\,,
      & \sigma \in \Sigma\,, \\
      \bot\,, & \sigma \not\in \Sigma\,.
    \end{cases}
  \end{align*}
  For $\sigma \in \Sigma$, observe that $\chi_{\sigma}(z) = \bot$ if
  $Y_{\sigma}$ is empty and otherwise that $\chi_{\sigma}(z) = \top$
  if and only if $\chi_{\sigma}^{y}(z) = \top$ for some $y \in
  Y_{\sigma}$.
  
  We define a modal algebra structure on the product Boolean algebra
  $\Alg{A} \times 2$. The modal operators $\pos[\sigma]$ are defined
  by:
  \begin{align*}
    \pos[\sigma](z,w) & =
    (\pos[\sigma] z, \chi_{\sigma}(z))\,.
  \end{align*}
  Since the functions $\chi_{\sigma}$ preserve joins, these modal
  operators are normal (i.e. they preserve finite joins).  Also,
  observe that the first projection $\prj_{1} : \Alg{A}\times 2 \rTo
  \Alg{A}$ is a morphism of modal algebras.
  
  Suppose now that $\Alg{A}$ is freely generated by a set $P$,
  $\Alg{A} = \Alg{F}_{P}$ , and let $\Lambda$ be a set of literals
  such that $\bigwedge \Lambda\not\leq \bot$.  Since $p$ and $\neg p$
  cannot belong both to $\Lambda$, we can choose a function $f :
  P \rTo \mathcal{A}\times 2$ with these properties: (i) $f(p) \in
  \set{(p,\bot),(p,\top)}$ for each $p \in P$, (ii) $f(p) = (p,\top)$
  if $p \in \Lambda$ and $f(p) = (p,\bot)$ if $\neg p \in \Lambda$.
  
  Let $ \tilde{f} : \Alg{F}_{P} \rTo \Alg{F}_{P}\times 2 $ be the
  extension of $f$ to a modal-algebra homomorphism, and observe that
  $\prj_{1} \circ \tilde{f} = \id_{\Alg{F}_{P}}$, since this relation
  holds on generators, and that $\tilde{f}(l) = (l,\top)$ for $l \in
  \Lambda$. Suppose that
  $$
  \bigwedge \Lambda
  \land
  \bigwedge_{\sigma \in \Sigma} 
  (\nec[\sigma]\bigvee Y_{\sigma}  \land
  \bigwedge_{y \in Y_{\sigma}}  \pos[\sigma] y)
  \leq \bot\,.
  $$
  If we
  apply the morphism $\tilde{f}$ to the above expression we obtain
  $$
  (\,
  \bigwedge \Lambda
  \land
  \bigwedge_{\sigma \in \Sigma} 
  (\nec[\sigma]\bigvee Y_{\sigma}  \land
  \bigwedge_{y \in Y_{\sigma}}  \pos[\sigma] y)
  \,
  ,\,
  a \land \bigwedge_{\sigma \in \Sigma}
  (b_{\sigma}\land c_{\sigma})\,)
  \leq (\bot,\bot)\,,
  $$
  where
  \begin{align*}
    a & = \bigwedge_{l \in \Lambda} \prj_{2}(\tilde{f}(l)) 
    = \bigwedge_{l \in \Lambda} \top \
    = \top\,,
    \tag*{since $\tilde{f}(l) = (l,\top)$}
    \\
    b_{\sigma}
    &= \neg \chi_{\sigma}(\neg \bigvee Y_{\sigma} ) 
    = \top\,,
    \intertext{%
      -- this relation is trivial if $Y_{\sigma}$ is empty, and
      otherwise note that $\chi_{\sigma}(\neg \bigvee Y_{\sigma})
      = \top$ iff $\chi^{y}_{\sigma}(\neg \bigvee Y_{\sigma}) =
      \top$ for some $y \in Y_{\sigma}$,
      which cannot be because of
      $
      \bot = \chi^{y}_{\sigma}(\bot) =
      \chi^{y}_{\sigma}(y \land \neg \bigvee
      Y_{\sigma}) = \chi^{y}_{\sigma}(y) \land 
      \chi^{y}_{\sigma}(\neg \bigvee
      Y_{\sigma}) = \top
      $ -- and finally 
    }
    c_{\sigma}
    & = \bigwedge_{y \in Y_{\sigma}} 
    \chi_{\sigma}(y)  =  \bigwedge_{y \in Y_{\sigma}} \top
    = \top\,. 
  \end{align*}
  We obtain
  $
  a \land \bigwedge_{\sigma \in \Sigma}
  b_{\sigma}\land c_{\sigma}
  =  \top 
  $ 
  which contradicts 
  $
  a \land \bigwedge_{\sigma \in \Sigma}
  b_{\sigma}\land c_{\sigma}
  \leq  \bot$.
  
\end{proof}
We extend now the previous result from modal algebras to modal
$\mu$-algebras.
\begin{Proposition}
  The implication \eqref{eq:mainresult} holds in a free modal
  $\mu$-algebra.
\end{Proposition}
\begin{proof}
  The proposition follows since if $\Alg{A}$ is a modal $\mu$-algebra,
  then the modal algebra $\Alg{A}\times 2$ is also a modal
  $\mu$-algebra and the first projection is a morphism of modal
  $\mu$-algebras.  This can be seen as follows: suppose that we have
  defined the interpretation of a term $f$ in the algebra
  $\Alg{A}\times 2$ as an operation $f : (\Alg{A}\times 2)^{\set{x}
    \cup Y} \rTo \Alg{A}\times 2$ so that the first projection
  preserves the interpretation. This is equivalent to saying that, for
  any fixed $v \in (\Alg{A}\times 2)^{Y}$, the following diagram
  commutes:
  $$
  \mydiagram[6em]{
    [](!s{\Alg{A}\times 2}{\Alg{A}}{\Alg{A}\times 2}{\Alg{A}}
    {1}{1},
    !a{^{\prj}}
    {^{f_{v}}}
    {^{f_{\prj(v)}}}
    {^{\prj}}
    )
  }
  $$
  Then $f_{v} = \langle \prj \circ f_{v},\psi \rangle = \langle
  f_{\prj(v)} \circ \prj,\psi \rangle$ for some $\psi : \Alg{A}\times
  2 \rTo 2$.  Considering that for each fixed $a \in \Alg{A}$
  $\mu_{y}.\psi(a,y)$ exists -- $2$ is a complete lattice -- we can
  use the Beki\v{c} property to argue that the least fixed point of
  $f_{v}$ exists and is equal to the pair
  $(\mu_{x}.f_{\prj(v)},\mu_{y}.\psi(\mu_{x}.f_{\prj(v)},y))$.
  Therefore we interpret the term $\mu_{x}.f$ in $\Alg{A}\times 2$ as
  suggested above, so that the first projection $\prj$ preserves the
  interpretation of the term $\mu_{x}.f$.
  
  Since all the terms of the theory of modal $\mu$-algebras are
  generated either by substitution or by formation of fixed points
  from the terms of the theory of multi-modal algebras, we deduce
  that $\Alg{A}\times 2$ is a modal $\mu$-algebra.  
\end{proof}

\begin{Lemma}
  \label{lemma:equivalence}
  On any modal algebra $\Alg{A}$ condition \eqref{eq:mainresult} 
  is equivalent to
  \begin{equation}
    \label{eq:smainresult}
    \begin{split}
      \bigwedge \Lambda \land \bigwedge_{\sigma \in Act}
      &(\nec[\sigma] x_{\sigma} \land\bigwedge_{y \in Y_{\sigma}}
      \pos[\sigma] y )
       \leq \bot \\
      &\timplies \\
      \bigwedge \Lambda &\leq \bot 
      \tor
      \exists \sigma \in Act, y \in
      Y_{\sigma} \tst  x_{\sigma}\land y\leq \bot\,.
    \end{split}
  \end{equation}
\end{Lemma}
\begin{proof}
  Assume that \eqref{eq:smainresult} holds and that the antecedent of
  \eqref{eq:mainresult} holds for some $\Sigma \subseteq Act$ and sets
  $Y_{\sigma}$.  In \eqref{eq:smainresult} let $x_{\sigma} =\bigvee
  Y_{\sigma}$ if $\sigma \in \Sigma$ and $x_{\sigma} = \top$ and
  $Y_{\sigma} = \emptyset$ if $\sigma \not\in \Sigma$. It immediately
  follows that either $\bigwedge \Lambda \leq \bot$, or there exists
  $\sigma \in \Sigma$ and some $y \in Y_{\sigma}$ such that $y \leq
  \bot$.
  
  Conversely, assume that
  $$
  \bigwedge \Lambda  \land
  \bigwedge_{\sigma \in Act} 
  (\nec[\sigma] x_{\sigma} \land
  \bigwedge_{y \in Y_{\sigma}} \pos[\sigma] y) 
  \leq \bot
  $$
  and derive
  \begin{align*}
    \bigwedge \Lambda \land
    \bigwedge_{\sigma \in Act} 
    (\,
    \nec[\sigma]((\bigvee Y_{\sigma}) \land x_{\sigma})
    \land
    \bigwedge_{y \in Y_{\sigma}} 
    \pos[\sigma](y \land x_{\sigma}) 
    \,)
    = \hspace{10mm}
    &  \\
    \bigwedge \Lambda  \land
    \bigwedge_{\sigma \in Act} 
    (\,
    \nec[\sigma](\bigvee_{y
      \in Y_{\sigma}} ( y \land x_{\sigma}) )
    \land
    \bigwedge_{y \in Y_{\sigma}}
    \pos[\sigma](y \land x_{\sigma}) \,)
    \leq & 
    \bot\,,
  \end{align*}
  using the fact that all operations involved are order
  preserving and distributivity.
  If we also assume that condition \eqref{eq:mainresult} holds, then it follows that $\bigwedge \Lambda \leq \bot$ or $x \land y_{\sigma}
  \leq \bot$ for some $\sigma \in Act$ and $y \in Y_{\sigma}$.
  
\end{proof}
The last Proposition has almost lead us to a proof of Theorem
\ref{theo:mainsec}.  In order to complete the proof, we need to argue
that if $\bigwedge \Lambda \leq \bot$ in a free modal $\mu$-algebra,
then $p,\neg p \in \Lambda$ for some generator $p$.  However the
latter property holds in a Boolean algebra $\mathcal{B}_{P}$ freely
generated by the set $P$, so that it is enough to argue that the
unique Boolean algebra homomorphism $\kappa : \mathcal{B}_{P} \rTo
\mathcal{F}_{P}$ extending the inclusion of generators $j : P \rTo
\mathcal{F}_{P}$ is an embedding. To this goal, observe that we can
assume $P$ to be finite so that the Boolean algebra $\mathcal{B}_{P}$
is finite as well, hence it is complete. $\mathcal{B}_{P}$ can also be
given a trivial structure of a modal algebra (say $\pos[\sigma]x = x$)
and therefore it is a modal $\mu$-algebra.  Let $\tilde{f}:
\mathcal{F}_{P} \rTo \mathcal{B}_{P}$ be the morphism of modal
$\mu$-algebras such that $\tilde{f}\circ j(p) = p$ for $p \in P$, then
$\tilde{f} \circ \kappa = \id_{\mathcal{B}_{P}}$, since this relation
holds on generators, and $\kappa$ is an embedding.

\section{First consequences}

In this section we present the first consequences of the property
stated in Theorem \ref{theo:mainproperty}.  We shall prove Claim
\ref{claim:uno} stating that modal operators $\pos[\sigma]$ are left
adjoints. Later we shall prove that a Kleene star modality,
$\pos[\sigma^{\ast}]$ in PDL notation, is constructive. This means
that this operation is a parametrized least prefixed point which is
the supremum over the chain of finite approximants. A proof of this
fact is included since it well exemplifies the theory that we shall
develop in the next sections.

\subsection*{Modal operators are adjoints}

Claim \ref{claim:uno} can also be understood by saying that reverse or
backward modalities are definable in free modal $\mu$-algebras.  This
property is analogous to Brzozowski derivatives being definable on
free Kleene-algebras \cite{kozen:MyN} and part of our contributition
consists in adapting the ideas presented there to the context of the
propositional modal $\mu$-calculus.
\begin{Proposition}[i.e. Claim \ref{claim:uno}]
  \label{prop:adj}
  On a free modal $\mu$-algebra each modal operator $\pos[\sigma]$ is
  a left adjoint.
\end{Proposition}
\begin{proof}
  Each element of a free modal $\mu$-algebra is a meet of elements of
  the form $\bigvee \Lambda \vee \bigvee_{\tau \in Act} (\,\pos[\tau]
  x_{\tau} \vee \bigvee_{y \in Y_{\tau}} \nec[\tau] y\,)$ where
  $\Lambda$ is a set of literals.  The previous statement holds since
  every term of the modal $\mu$-calculus is provably equivalent to a
  guarded term, see \cite{kozen}, i.e. to a term where negation
  appears only in front of generators and every bound fixed point
  variable is in the scope of some modal operator.  Using fixed point
  equalities it is possible to unravel the term to extract its first
  modal level.  The statement then follows by distributivity.
  
  Therefore, we begin by defining the right adjoint for an element
  having this form: if
  \begin{align*}
    b & = \bigvee \Lambda \vee \bigvee_{\tau
    \in Act} (\,\pos[\tau] x_{\tau} \vee \bigvee_{y \in Y_{\tau}}
  \nec[\tau] y\,)\,,
  \end{align*} 
  then we define
  \begin{align*}
    \righta[\sigma](b) & =
    \begin{cases}
      \top\,, & \text{if } b = \top\,,\\
      x_{\sigma}\,, & \text{otherwise}.
    \end{cases}
  \end{align*}
  We argue now that $\pos[\sigma]x \leq b$ iff $x
  \leq\righta[\sigma](b)$. Suppose that $\pos[\sigma]x \leq b$: if $b
  = \top$ then clearly $x \leq \top = \righta[\sigma](x)$, and if $b
  \neq \top$, then we deduce $x \leq x_{\sigma} = \righta[\sigma](b)$.
  The latter statement is a consequence of Theorem
  \ref{theo:mainproperty} when properly dualized, taking into account
  that all the disjuncts other than $x \land \neg x_{\sigma} \leq
  \bot$ in the consequent of \ref{theo:mainproperty} imply $b = \top$.
  Conversely, the relation $\pos[\sigma]\righta[\sigma](b) \leq b$
  clearly holds and implies that $x
  \leq \righta[\sigma](b)$ implies $\pos[\sigma]x \leq b$.  Note
  also that $\righta[\sigma](b)$ does not depend on the representation
  of $b$, as it is uniquely determined by the property $x \leq
  \righta[\sigma](b)$ iff $\pos[\sigma]x \leq b$.
  
  It is a standard step then to extend the right adjoint to all the
  elements of a free modal $\mu$-algebra: if $x = \bigwedge_{j \in J}
  b_{j}$, then we define $\righta[\sigma](x) = \bigwedge_{j \in
    J}\righta[\sigma](b_{j})$.  
\end{proof}

\subsection*{The Kleene star is constructive}

An important property of $\righta[\sigma](z)$ -- the right adjoint to
$\pos[\sigma]$ defined in the proof of Proposition \ref{prop:adj} -- is
that it is computed out of the syntax of $z$. More precisely,
$\righta[\sigma](z)$ is computed as a meet of terms belonging to the
Fisher-Ladner closure, see \cite{kozen}, of a term representing $z$.
The Fisher-Ladner closure has to be thought as the space of subterms
of $z$, in particular it is finite.  Consequently, the set
$\set{\righta[\sigma]^{n}(z) \,|\, n \geq 0}$ is finite and
$\bigwedge_{n \geq 0} \righta[\sigma]^{n}(z)$ exists in a free modal
$\mu$-algebra. We exemplify how to exploit this fact by proving that
$\mu_{y}.(x \vee \pos[\sigma]y)$ is the supremum over the chain of its
finite approximants. 

We shall use the standard Propositional Dynamic
Logic notation and let $\pos[\sigma^{\ast}]x = \mu_{y}.(x \vee
\pos[\sigma^{\ast}]y)$.
\begin{Lemma}
  The relation
  \begin{align*}
    \pos[\sigma^{\ast}]a & = \bigvee_{n \geq 0} \pos[\sigma]^{n}a
  \end{align*}
  holds in a free modal $\mu$-algebra.
\end{Lemma}
\begin{proof}
  We only need to prove that if $\pos[\sigma]^{n}a
  \leq b$ for each $n \geq 0$, then $\pos[\sigma^{\ast}]a \leq b$. 

  Assume that $\pos[\sigma]^{n}a
  \leq b$ for each $n \geq 0$ and
  transpose these relations to obtain $a \leq
  \righta[\sigma]^{n}(b)$ for each $n \geq 0$, hence $a \leq
  \bigwedge_{n\geq 0} \righta[\sigma]^{n}(b)$.
  We claim that $\bigwedge_{n\geq 0} \righta[\sigma]^{n}(b)$ is a
  $\pos[\sigma]$-prefixed point. Indeed:
  \begin{align*}
    \pos[\sigma] \bigwedge_{n\geq 0} \righta[\sigma]^{n}(b) & \leq
    \bigwedge_{n\geq 0} \pos[\sigma] \righta[\sigma]^{n}(b)
    \tag*{$\pos[\sigma]$ is order preserving}
    \\
    & = \pos[\sigma]b \land \bigwedge_{n\geq 0} \pos[\sigma]
    \righta[\sigma]^{n+1}(b)
    \\
    & \leq \pos[\sigma]b \land \bigwedge_{n\geq 0}
    \righta[\sigma]^{n}(b) \leq \bigwedge_{n\geq 0}
    \righta[\sigma]^{n}(b) \tag*{by the counit relation
      $\pos[\sigma]\righta[\sigma]x \leq x$.}
  \end{align*}
  Thus $\bigwedge_{n\geq 0} \righta[\sigma]^{n}(b)$ is a
  $\pos[\sigma]$-prefixed point above $a$ and therefore
  $\pos[\sigma^{\ast}] a \leq \bigwedge_{n\geq 0} \righta[\sigma]^{n}
  (b)$. Since $\bigwedge_{n\geq 0}\righta[\sigma]^{n}(b) \leq b$ we
  deduce $\pos[\sigma^{\ast}] a \leq b$. 
\end{proof}


\section{$\Opens$-adjoints of finite type}

The proof that the parametrized least prefixed point corresponding to
the PDL star modality $\pos[\sigma^{\ast}]$ is the supremum over the
chain of its finite approximants relies on the modality $\pos[\sigma]$
being a left adjoint. We cannot use this idea on the nose to prove
constructiveness of other operations that are not left adjoints.  For
example, a necessity modal operation $\nec[\sigma]$ is not a left
adjoint on free modal $\mu$-algebras since it doesn't preserve joins.
To deal with the general case left adjoints are generalized as
follows.
\begin{Definition}
  Let $L$ and $M$ be posets.  An order preserving function $f : L \rTo
  M$ is a left \emph{$\Opens$-adjoint} if for each $m \in M$ the set
  \begin{align*}
    \set{x\,|\, f(x) \leq m\,}
  \end{align*}
  is a finitely generated lower set. 
\end{Definition}
That is, $f$ is a $\Opens$-adjoint iff the above set is a finite union
of principal ideals, or equivalently iff for each $m \in M$ there
exists a finite set $\Cvrs{f}{m}$ such that for all $x \in L$ $f(x)
\leq m$ if and only if $x \leq c$ for some $c \in \Cvrs{f}{m}$.  We
shall say that $\Cvrs{f}{m}$ is the set of $f$-covers of $m$ or the
covering set of $f$ and $m$.

It is easily seen that $f$ is a left adjoint if and only if
$\set{x\,|\, f(x) \leq m\,}$ is a principal ideal, thus every left
adjoint is a left $\Opens$-adjoint. Also,  $f$
is a left $\Opens$-adjoint if and only if
$$
\Opens(f): \Opens(L) \rTo \Opens(M)
$$
is a left adjoint; here $ \Opens(P)$ is the set of finitely
generated lower sets of the poset $P$ and $\Opens(f)$ is the obvious
map induced by this functorial construction. The notion of
$\Opens$-adjoint presented here corresponds to that of
a $Pro(\mathcal{D})$-adjoint \cite{tholen} where $\mathcal{D}$ is the
class of all finite discrete categories.  Similar but slightly
different is the notion of a multiadjoint \cite{diers}. In the
following, $\Opens$-adjoint will abbreviate left $\Opens$-adjoint.

We begin presenting an interesting order theoretic property of
$\Opens$-adjoints:
\begin{Lemma}
  A $\Opens$-adjoint $f$ is continuous: if $I$ is a directed set and
  $\bigvee I$ exists, then $\bigvee_{i \in I} f(i)$ exists as well and
  is equal to $f(\bigvee I)$.
\end{Lemma}
\begin{proof}
  Suppose that for all $i \in I$ $f(i) \leq m$. We can find $c_{i} \in
  \Cvrs{f}{m}$ such that $i \leq c_{i}$.  Since $I$ is directed and
  the $c_{i}$ are finite, we can find $i_{0}$ such that $i \leq
  c_{i_{0}}$ for all $i \in I$ and consequently $\bigvee I \leq
  c_{i_{0}}$. It follows that $f(\bigvee I) \leq f(c_{i_{0}}) \leq m$.
  
\end{proof}
We can argue that being a $\Opens$-adjoint is a stronger property than
merely being continuous by considering the binary meet $\land :
\mathcal{B} \times \mathcal{B} \rTo \mathcal{B}$ on an infinite
Boolean algebra $\mathcal{B}$. The binary meet is continuous -- since
it is continuous in each variable -- but it is not a $\Opens$-adjoint.
This can be seen by computing a candidate covering set
$\Cvrs{\land}{\bot}$. Since $x \land \neg x \leq \bot$, then we should
be able to find $(\alpha_{x},\beta_{x}) \in \Cvrs{f}{\bot}$ such that
$x \leq \alpha_{x}$, $\neg x \leq \beta_{x}$, and moreover $\alpha_{x}
\land \beta_{x} \leq \bot$. It follows that $\alpha_{x} \leq \neg
\beta_{x} \leq x$ and $\alpha_{x} = x$. Thus, for an infinite Boolean
algebra the covering set $\Cvrs{\land}{\bot}$ has to be infinite.

We list next some properties of $\Opens$-adjoints:
\begin{Proposition}\mbox{\hspace{0mm}}
  \label{prop:list}
  \begin{enumerate}
  \item \label{item:leftadj} An order preserving function $f: L \rTo
    M$ is a left adjoint if and only if it is a $\Opens$-adjoint and
    preserves finite joins.
  \item \label{item:meetgen} If a lattice $M$ is finitely
    meet-generated by a subset $B \subseteq M$, then $f: L \rTo M$ is
    a $\Opens$-adjoint if and only if the covering set $\Cvrs{f}{b}$
    exists for each $b \in B$.
  \item \label{item:comps} Identities are $\Opens$-adjoints, and
    $\Opens$-adjoints are closed under composition.
  \item \label{item:prods} If the domain posets are meet semilattices,
    then the projections $\prj_{i} : L_{1} \times L_{2} \rTo L_{i}$,
    $i = 1 ,2$, are $\Opens$-adjoints. Moreover $\langle
    f_{1},f_{2}\rangle : L \rTo M_{1}\times M_{2}$ is a
    $\Opens$-adjoint provided that $f_{i} : L \rTo M_{i}$, $i = 1 ,2$,
    are $\Opens$-adjoints.
  \item \label{item:alladjs} Finite joins are $\Opens$-adjoints. 
  \item \label{item:const} Constant functions are $\Opens$-adjoints.
    If $L$ is an Heyting algebra (or a Browerian semilattice), then
    $f(x) = k \land x : L \rTo L$ is a $\Opens$-adjoint, where $k$ is
    a constant.
  \end{enumerate}
\end{Proposition}
\begin{proof}
  \ref{item:leftadj}: If $\righta[f]$ is right adjoint to $f$, then
  the lower set $\set{y\,|\,f(y) \leq m}$ is generated by
  $\righta[f](m)$, thus $f$ is a $\Opens$-adjoint.  For the second
  statement, define the right adjoint $\righta[f](m)$ as $\bigvee
  \Cvrs{f}{m}$.  

  \ref{item:meetgen}: Let $m = \bigwedge_{i \in I} b_{i}$. If $f(x)
  \leq m$, then $f(x) \leq b_{i}$ for all $i \in I$ and there exists
  $c_{i} \in \Cvrs{f}{b_{i}}$ such that $x \leq c_{i}$: therefore $x
  \leq \bigwedge_{i \in I} c_{i}$.  Conversely, if $x \leq
  \bigwedge_{i \in I} c_{i}$ with $c_{i} \in \Cvrs{f}{b_{i}}$ for each
  $i \in I$, then $f(x) \leq b_{i}$, $i \in I$, and $f(x) \leq m$.
  That is, we can define
  \begin{align*}
    \Cvrs{f}{\bigwedge_{i \in I} b_{i}}
    & = \bigwedge_{i \in I} \Cvrs{f}{b_{i}}\,.
  \end{align*}
  This set is finite if $I$ is finite.  

  \ref{item:comps}:
  The identity is left adjoint to itself.  For composition we can
  define:
  \begin{align*}
    \Cvrs{f \circ g}{m}
    & = \bigcup_{ c \in \Cvrs{f}{m}} \Cvrs{g}{c} \,.
  \end{align*}

  \ref{item:prods}:
  Since we are assuming existence of $\top$, projection functions  are
  left adjoints.  For pairing we define:
  \begin{align*}
    \Cvrs{\langle f_{1},f_{2}\rangle}{(m_{1},m_{2})}
    & = \set{ c_{1} \land c_{2} \,|\,
      c_{1} \in  \Cvrs{f_{1}}{m_{1}},\,
      c_{2} \in \Cvrs{f_{2}}{m_{2}} }\,.
  \end{align*}
  
  \ref{item:alladjs}: The diagonal is right adjoint to $\vee : L
  \times L \rTo L$.

  \ref{item:const}:
  Let $f_{k}$ be the constant function taking every $x$ to the
  constant value $k$. We can define
  \begin{align*}
    \Cvrs{f_{k}}{m} & = 
    \begin{cases}
      \emptyset & k \not\leq m \\
      \set{\top} & \textrm{otherwise}.
    \end{cases}
  \end{align*}
  The operation $k \rightarrow y$ is right adjoint to $k \land
  x$.  
\end{proof}

\subsection*{$\Opens$-adjoints and fixed points}

We analyze next $\Opens$-adjoints for which it makes sense to consider
least fixed points, i.e. those of the form $f: L^{x}\times M^{y} \rTo
L$. For such an $f$, we define a directed multi-graph
$\mathcal{G}_{x}(f,L)$ as follows:
\begin{itemize}
\item its vertices are elements of $L$,
\item there is a transition $l \rTo^{m'} l'$ iff $(l',m') \in
  \Cvrs{f}{l}$.
\end{itemize}
We write $\mathcal{G}_{x}(f,l)$ for the full subgraph of
$\mathcal{G}_{x}(f,L)$ of elements of $L$ that are reachable from $l$:
$l' \in L$ is a vertex of $\mathcal{G}_{x}(f,l)$ iff there exists a
path from $l$ to $l'$ in $\mathcal{G}_{x}(f,L)$.
\begin{Definition}
  We say that the $\Opens$-adjoint $f: L^{x}\times M^{y} \rTo L$ has
  \emph{finite type} for the variable $x$ if for each $l \in L$ the
  graph $\mathcal{G}_{x}(f,l)$ is finite.
\end{Definition}

\begin{Lemma}
  Suppose that $M$ is a meet semilattice, the $\Opens$-adjoint $f:
  L^{x}\times M^{y} \rTo L$ has finite type, and $\mu_{x}.f(x,y)$
  exists for each $y \in M$. Then the order preserving parametrized
  fixed point $\mu_{x}.f: M^{y}
  \rTo L$ is again a $\Opens$-adjoint. 
\end{Lemma}
\begin{proof}
  Recall that a path of length $n$ in $\mathcal{G}_{x}(l,M)$ is a
  sequence of transitions $l_{i} \rTo^{m_{i + 1}} l_{i+1}$ with $0
  \leq i < n$.  Such a path is infinite if $n = \omega$. The path is
  from $l$ if $l_{0} = l$.
  
  Remark that in an infinite path $l_{i} \rTo^{m_{i + 1}} l_{i+1}$, $i
  < \omega$, there exists only a finite number of $m$'s such that $m =
  m_{i}$ for some $i$. Hence the meet $\bigwedge_{i \geq 1} m_{i}$
  exists in $M$.  We define
  \begin{align*}
    m \in \Cvrs{\mu_{x}.f}{l} & \text{ iff } m
    = \bigwedge_{i \geq 1} m_{i} \\
    & \hspace{1cm}\text{ for some infinite path }
    \set{l_{i} \rTo^{m_{i + 1}} l_{i +1}}_{i \geq 0} \text{ from } l\,.
  \end{align*}
  Observe that this set is actually finite, as a consequence of
  $\mathcal{G}_{x}(f,l)$ being finite. 
  
  We begin verifying that $\mu_{x}.f(m) \leq l$ if $m \in
  \Cvrs{\mu_{x}.f}{l}$.  Observe that, by monotonicity, $f(l_{i+1},m)
  \leq f(l_{i+1},m_{i + 1}) \leq l_{i}$ for all $i \geq 0$, and more
  generally $f_{m}^{k}(l_{i + k}) \leq l_{i}$ for all $i,k\geq 0$.
  Choose $i < j$ such that $l_{i} = l_{j}$ and let $k = j - i$, then $
  f^{k}_{m}(l_{i}) = f_{m}^{k}(l_{j}) \leq l_{i}$, hence
  $\mu_{x}.f(x,m) = \mu_{x}.f^{k}_{m}(x)\leq l_{i}$.  We deduce
  $\mu_{x}.f(x,m) = f^{i}_{m}(\mu_{x}.f(x,m)) \leq f^{i}_{m}(l_{i})
  \leq l_{0} = l$.
  
  Conversely, assume that $\mu_{x}.f(x,y) \leq l_{0}$: we can use the
  fixed point equation to deduce $f(\mu_{x}.f(x,y),y) \leq l_{0}$
  which in turn implies $(\mu_{x}.f(x,y),y) \leq (l_{1},m_{1})$ for
  some pair $(l_{1},m_{1}) \in \Cvrs{f}{l_{0}}$. By iterating the
  procedure, we can construct an infinite path $\set{l_{i} \rTo^{m_{i +
  1}} l_{i +1}}_{i \geq 0}$
  from $l$ such that for all $i \geq 1$ we have $(\mu_{x}.f(x,y) , y)
  \leq (l_{i},m_{i}) $. We have  therefore $y \leq \bigwedge_{i \geq 1}
  m_{i} \in \Cvrs{\mu_{x}.f}{l}$.  
\end{proof}
It is a natural step to prune covering sets $\Cvrs{f}{m}$ to extract
the antichain of maximal elements. If this operation is performed on
$\Cvrs{\mu_{x}.f}{l}$, we see that a maximal element is a meet indexed
by some pan in $\mathcal{G}_{x}(f,l)$. By a pan, we mean a finite path
that can be split into a simple path followed by a simple cycle.

\begin{Lemma}
  \label{lemma:finitaryconstructive}
  Under the conditions of the previous Lemma, the least prefixed
  point of $f : L^{x} \times M^{y} \rTo L$ is constructive:
  \begin{align*}
    \mu_{x}.f(x,y)
    & = \bigvee_{n\geq 0} f^{n}_{y}(\bot)\,.
  \end{align*}
\end{Lemma}
\begin{proof}
  Assume $l$ is such that $f^{n}_{y}(\bot) \leq l$ for each $n\geq0$.
  Let $k$ be the number of vertices in the graph
  $\mathcal{G}_{x}(f,l)$ and observe that the relation
  $f^{k}_{y}(\bot) \leq l$ implies that we can find a path of length
  $k$ $l_{i}\rTo^{m_{i + 1}} l_{i + 1}$ from $l$ with the property
  that $y \leq m_{i}$ for $i = 1,\ldots ,k$.  By choosing $i,j$ such
  that $0 \leq i < j \leq k$ and $l_{i} = l_{j}$, construct an
  infinite path $l_{i}\rTo^{m_{i + 1}} l_{i + 1}$ from $l$ such that
  $y \leq m_{i}$ for $i \geq 1$.
  
  Thus $\bigwedge m_{i} \in \Cvrs{\mu_{x}.f}{l}$ and
  therefore $\mu_{x}.f(x,y) \leq \mu_{x}.f(x,\bigwedge m_{i}) \leq l$.
  
\end{proof}

\subsection*{$\Opens$-adjoints on free modal $\mu$-algebras}

We continue by considering $\Opens$-adjoints on free modal
$\mu$-algebras.  We have seen that meets provide a counter-example for
$\Opens$-adjointness. In \cite{janin} the authors suggest a sort of
best approximation of meets as $\Opens$-adjoints.  They define the
arrow term by:
\begin{align}
  \label{eq:arrow}
  \arrow[\sigma] X & = \nec[\sigma]\bigvee X \land \bigwedge_{x \in X}
  \pos[\sigma] x\,,
\end{align}
and, for a set of literals $\Lambda$, for a subset $\Sigma \subseteq
Act$, and for disjoint sets of variables $\set{X_{\sigma}}_{\sigma \in
  \Sigma}$, they also define the special conjunction term by:
\begin{align}
  \label{eq:spcon}
  \spcon[\Lambda,\Sigma] \set{X_{\sigma}} & = \Lambda \land
  \bigwedge_{\sigma \in \Sigma} \arrow[\sigma] X_{\sigma}\,.
\end{align}
Let $X = \bigcup_{\sigma \in \Sigma} X_{\sigma}$ and $v \in
\mathcal{F}^{X}$ be a vector of elements of a free modal
$\mu$-algebra. We have seen in \ref{lemma:equivalence} that
$\spcon[\Lambda,\Sigma] v = \bot$ if either the literals in $\Lambda$
are inconsistent or $v(x) = \bot$ for some $\sigma \in \Sigma$ and $x
\in X_{\sigma}$.
\begin{Lemma}
  Special conjunctions on free modal $\mu$-algebras are
  $\Opens$-adjoints of finite type.
  \label{lemma:specconj}
\end{Lemma}
\begin{proof}
  Recall from \ref{prop:adj} that the free modal $\mu$-algebra is
  finitely meet-generated by elements of the form 
  \begin{align}
    \label{eq:specialels}
    b & = \bigvee \Gamma \vee \bigvee_{\tau
    \in Act}( \pos[\tau] d_{\tau} \vee \bigvee_{e \in
    E_{\tau}} \nec[\tau] e)\,,
  \end{align}
  where $\Gamma$ is a set of literals.  By Proposition
  \ref{prop:list}.\ref{item:meetgen}, it is enough to define the
  covering sets $\Cvrs{\spcon[\Lambda,\Sigma]}{b}$ for such $b$'s.
  Observe that, if $b = \top$, then we can define $\Cvrs{f}{\top} =
  \set{\top}$ for any monotone $f$.  Also, if $\bigwedge \Lambda \leq
  \bigvee \Gamma$, then we can define
  $\Cvrs{\spcon[\Lambda,\Sigma]}{b} = \set{\top}$.
  
  Hence, let $b$ be as in \eqref{eq:specialels} and suppose that $b
  \neq \top$ and $\bigwedge \Lambda \not\leq \bigvee \Gamma$.
  Recalling that $X$ is the disjoint union of the $X_{\sigma}$,
  $\sigma \in \Sigma$, we define
  \begin{align*}
    \Cvrs{\spcon[\Lambda,\Sigma]}{b} & = 
    \set{c_{\sigma,y} \,|\, \sigma \in \Sigma, y \in X_{\sigma}} 
    \cup \set{
      c_{\sigma,e} \,| \,\sigma \in \Sigma, e \in E_{\sigma}}\,,
  \end{align*}
  where the vectors $c_{\sigma,y}, c_{\sigma,e} \in \mathcal{F}^{X}$
  are as follows:
  \begin{align*}
    c_{\sigma,y}(x)
    & = 
    \begin{cases}
      d_{\sigma}, & x = y\,,\\
      \top, & \text{otherwise},      
    \end{cases}
    \intertext{and}
    c_{\sigma,e}(x)
    & = 
    \begin{cases}
      \top, & x \in X_{\tau}, \tau \neq \sigma \\
      d_{\sigma} \vee e, & x \in X_{\sigma}\,.
    \end{cases}
  \end{align*}
  Observe that 
  \begin{align*}
    \spcon[\Lambda,\Sigma](c_{\sigma,y})
    & \leq 
    \pos[\sigma]d_{\sigma} \leq b\,,
    \intertext{and}
    \spcon[\Lambda,\Sigma](c_{\sigma,e})
    &\leq\; \arrow[\sigma]\set{d_{\sigma} \vee e} 
    \;\leq\; \nec[\sigma](d_{\sigma} \vee e) 
    \; \leq\; \pos[\sigma]d_{\sigma} \vee \nec[\sigma]e 
    \;\leq \;b\,.
  \end{align*}
  It follows that if $v \leq c \in
  \Cvrs{\spcon[\Lambda,\Sigma]}{b}$, then $\spcon[\Lambda,\Sigma] v
  \leq b$.
  
  Conversely, let $v \in \mathcal{F}^{X}$, and suppose that
  $\spcon[\Lambda,\Sigma](v) \leq b$. We apply Theorem
  \ref{theo:mainproperty} to this relation, whose explicit expression
  is
  $$
  \bigwedge \Lambda \;\land\; \bigwedge_{\sigma \in \Sigma}
  (\nec[\sigma]\bigvee_{x \in X_{\sigma}} v(x)  
  \land \bigwedge_{x \in X_{\sigma}}
  \pos[\sigma] v(x))
  \; \leq \;\bigvee \Gamma \vee \bigvee_{\tau
    \in Act}( \pos[\tau] d_{\tau} \vee \bigvee_{e \in
    E_{\tau}} \nec[\tau] e)\,.
  $$  
  Since $b \neq \top$ and $\bigwedge \Lambda \not\leq \bigvee \Gamma$,
  one of the following two cases holds:
  \begin{enumerate}
  \item there exists $\sigma \in \Sigma$ and $x \in X_{\sigma}$ such
    that $v(x) \leq d_{\sigma}$: in this case $v \leq c_{\sigma,x}$,
  \item there exists $\sigma \in \Sigma$ and $e \in E_{\sigma}$ such
    that $v(x) \leq d_{\sigma} \vee e$ for each $x \in X_{\sigma}$, in
    this case $v \leq c_{\sigma,e}$.
  \end{enumerate}

  To end the proof, we remark that covers of an element $c \in
  \mathcal{F}$ are meets of subterms of a term representing $c$,
  showing that special conjunctions have finite type.  
\end{proof}
It is now easy to argue that the necessity modal operation
$\nec[\sigma]$ is a $\Opens$-adjoint on a free modal $\mu$-algebra.
By Proposition \ref{prop:list}, this is a consequence of
$\nec[\sigma]$ belonging to the cone generated by joins and special
conjunctions, since the relation $\nec[\sigma]x =
\arrow[\sigma]\set{x} \vee \arrow[\sigma]\emptyset$ holds on every
modal algebra.


\subsection*{Uniform families of  $\Opens$-adjoint of finite type}
In the previous subsection we have studied properties of
$\Opens$-adjoint and seen that having finite type is quite relevant
for least fixed points. We develop next some tools by which it will be
easier to compute the type of a $\Opens$-adjoint.

A function scheme is a triple $(f,X,Y)$: intuitively $f$ is a function
symbol, $X$ and $Y$ are finite sets of variables, $X$ being the arity
of $f$ and $Y$ being its coarity. That is, $f$ is meant to represent a
function of the form $f : L^{X} \rTo L^{Y}$. For a function scheme
$(f,X,Y)$, an $f$-automaton is a pair $\langle Q,\Delta_{f} \rangle$
where $Q$ is a set of states and $\Delta_{f} \subseteq Q^{Y} \times
Q^{X}$. For a family $\mathcal{F}$ of function schemes an
$\mathcal{F}$-automaton is a tuple $\langle Q,\set{\Delta_{f}}_{f \in
  \mathcal{F}} \rangle$, where $\langle Q,\Delta_{f}
\rangle$ is an $f$-automaton for each $f \in \mathcal{F}$.  Let
$\langle Q,\set{\Delta_{f}}_{f \in \mathcal{F}} \rangle$ be an
$\mathcal{F}$-automaton and $Q_{0} \subseteq Q$: we let
$\reach{\langle Q,\set{\Delta_{f}}\rangle ,Q_{0}} = \langle P,\set{\Delta'_{f}}
\rangle$ be the least sub-$\mathcal{F}$-automaton of $\langle Q,
\set{\Delta_{f}} \rangle$ such that $Q_{0} \subseteq P$ and $v \in
P^{Y}$ and $v\Delta_{f} w $ implies $w \in P^{X}$ for each function
scheme $(f,X,Y) \in \mathcal{F}$. The relations $\Delta_{f}'$ are the
restriction of the $\Delta_{f}$ to $P$.

For a family $\mathcal{F}$ of $\Opens$-adjoints of the form $f:
L^{X}\rTo L^{Y}$, the $\mathcal{F}$-automaton
$\mathcal{A}_{\mathcal{F}}$ is defined as follows: its set of states
is $L$ and $v \Delta_{f} c$ iff $c \in \Cvrs{f}{v}$.
\begin{Definition}
  A family $\mathcal{F}$ of $\Opens$-adjoints is a \emph{uniform
    family of finite type} if the underlying set of the
  $\mathcal{F}$-automaton $\reach{\mathcal{A}_{\mathcal{F}},Q_{0}}$ is
  finite whenever $Q_{0} \subseteq L$ is finite.
\end{Definition}
The obvious reason to introduce this notion is:
\begin{Lemma}
  \label{lemma:uniformfinitary}
  Let $\mathcal{F}$ be a uniform family of finite type. If $f \in
  \mathcal{F}$ and $f: L^{X}\times L^{Y}\rTo L^{X}$, then $f$ has
  finite type for the variable $X$.
\end{Lemma}
\begin{proof}
  Let $v_{0} \in L^{X}$ and $S$ be the set of states of
  $\mathcal{G}_{X}(f,v_{0})$.  If $Q_{0} = \set{ l\,|\, v_{0}(x) = l
    \text{ for some } x \in X }$, then we claim that $S \subseteq
  P^{X}$, where $P$ is the underlying set of
  $\reach{\mathcal{A}_{\mathcal{F}},Q_{0}}$.  Indeed, $v_{0} \in
  P^{X}$ and if $v \in P^{X}$ and $v \rTo^{m'} v'$, then $v
  \Delta_{f}(v',m')$ so that $(v',m') \in P^{X}\times P^{Y}$ and $v'
  \in P^{X}$.  
\end{proof}

We investigate now closure properties of a family $\mathcal{F}$ w.r.t.
finiteness.
\begin{Lemma}
  If $\mathcal{F}$ is a uniform family of finite type and $x \in X$,
  then $\mathcal{F} \cup \set{\prj[X]_{x}}$ is a uniform family of
  finite type.
\end{Lemma}
\begin{proof}
  Recall that $\Cvrs{\prj[X]_{x}}{b} = \set{c}$ where $c(x) = b$ and
  $c(y) = \top$ for $y \neq x$.
  It is easily argued that the 
  underlying set of $\reach{\mathcal{A}_{\mathcal{F} \cup
      \set{\prj[X]_{x}}},Q_{0}}$ is contained in the underlying
  set of $\reach{\mathcal{A}_{\mathcal{F}},Q_{0}\cup \set{\top}}$.
  
\end{proof}
\begin{Lemma}
  \label{lemma:fincomp}
  If $\mathcal{F}$ is a uniform  family of finite type and $f,g \in
  \mathcal{F}$ with $f : L^{X} \rTo L^{Y}$ and $g: L^{Y} \rTo L^{Z}$,
  then $\mathcal{F} \cup \set{f\circ g}$ is a uniform family of finite
  type.
\end{Lemma}
\begin{proof}
  Recall that $\Cvrs{g\circ f}{b} = \bigcup_{c \in \Cvrs{g}{b}}
  \Cvrs{g\circ f}{c}$, from which it results that the 
  underlying set of $\reach{\mathcal{A}_{\mathcal{F} \cup
      \set{f \circ g}},Q_{0}}$ is the same as the underlying
  set of $\reach{\mathcal{A}_{\mathcal{F}},Q_{0}}$.
  
\end{proof}
By the previous Lemmas, we can always assume that a uniform family of
finite type $\mathcal{F}$ is closed under post-composition with
projections, that is, if $\langle f_{y}\rangle_{y \in Y} : L^{X} \rTo
L^{Y}$, then $f_{y} \in \mathcal{F}$ for each $y \in Y$.


\begin{Lemma}
  \label{lemma:uniformfamily}
  Let $\mathcal{F}$ be a uniform family of finite type which is closed
  under post-composition with projections, and let $f_{y} : L^{X} \rTo
  L$, $y \in Y$, be elements of $\mathcal{F}$.  Then $\mathcal{F} \cup
  \set{\langle f_{y}\rangle_{y \in Y}}$ is a uniform family of finite
  type.
\end{Lemma}
\begin{proof}
  Let $Q_{0} \subseteq L$ be a finite subset of $L$, and let 
  $\reach{\mathcal{A}_{\mathcal{F}},Q_{0}} = \langle
  Q,\set{\Delta_{f}} \rangle$, so that, by assumption, $Q$ is finite.
  
  Let $S$ be the meet-semilattice generated by $Q$ and let $P$ be the
  underlying set of $\reach{\mathcal{A}_{\mathcal{F} \cup \set{\langle
        f_{y}\rangle}}, Q_{0}}$: we claim that $P \subseteq S$.
  Clearly, $Q_{0} \subseteq Q \subseteq S$.
  
  We show now that for each $f: L^{X} \rTo L^{Y}$ in $\mathcal{F} \cup
  \set{\langle f_{y}\rangle}$, $w \in S^{Y}$ and $w \Delta_{f} c$
  implies $c \in S^{Y}$.
  
  We analyze first the case of a function of the form $f: L^{X} \rTo
  L$.  Let $w \in S$ and suppose that $w \Delta_{f} c$, i.e. $c \in
  \Cvrs{f}{w}$.  Since $w \in S$, we can write $w = \bigwedge w_{i}$
  where $w_{i} \in Q$. Hence $c = \bigwedge c_{i}$ where $c_{i} \in
  \Cvrs{f}{w_{i}}$: we have, therefore, $w_{i}\Delta_{f} c_{i}$ and
  $c_{i} \in Q^{X}$. Hence, $c = \bigwedge_{i} c_{i}$ belongs to
  $S^{X}$.
  
  We analyze now the case of a function of the form $f: L^{X} \rTo
  L^{Y}$, with $Y$ not a singleton. Thus $f = \langle f_{y} \rangle$
  where each $f_{y} : L^{X} \rTo L$ belongs to $\mathcal{F}$.  Let $w
  \in S^{Y}$, and suppose that $w \Delta_{\langle f_{y}\rangle} c$.
  This means that $c \in \Cvrs{\langle f_{y}\rangle}{w}$ so that $c =
  \bigwedge c_{y}$ where $c_{y} \in \Cvrs{f_{y}}{w(y)}$ for each $y
  \in Y$.  We have already argued that $c_{y} \in S^{X}$, hence 
  $c =  \bigwedge c_{y} \in S^{X}$ as well.  
\end{proof}


\section{Some constructive systems of equations}

An order preserving $F : \Alg{A}^{X}\times \Alg{A}^{Y} \rTo
\Alg{A}^{X}$ 
can be thought to be a system of equations whose least solution is
given by the least fixed point. The set $X = \set{x_{1},\ldots
  ,x_{n}}$ is the set of bound variables of the system and $Y =
\set{y_{1},\ldots ,y_{m}}$ is the set of free variables, the sets $X$
and $Y$ being disjoint. If $F = \langle F_{x} \rangle_{x \in X}$, then
we  represent such systems as expected:
$$
\begin{system}
  &\vdots \\
  x_{i} & = & F_{x_{i}}(x_{1},\ldots ,x_{n},y_{1},\ldots ,y_{m}) \\
  &\vdots
\end{system}
$$
The Beki\v{c} property ensures that such a system of equations has a
least solution in every modal $\mu$-algebra if each $F_{x}$ is the
interpretation of a term of the theory of modal $\mu$-algebras.

In this section we shall prove that, for many such $F$ on a free modal
$\mu$-algebra $\mathcal{F}$, the least prefixed point is the supremum
over the chain of its finite approximants. The results of the previous
sections allow us to easily derive this property for a restricted set
of systems called here disjunctive-simple. Then, we freely use ideas
and tools from \cite[\S 9]{AN} to enlarge the class of systems that
can be proved to be constructive.  An improvement w.r.t. this
monograph consists in adapting these tools in order to argue about
existence of infinite suprema and approximants.  Our last effort will
be to prove that all the systems $F = \langle F_{x} \rangle$ whose
$F_{x}$ are elementary operations of the theory of modal algebras
enjoy this property.



\begin{Definition}
  We say that a term of the theory of modal $\mu$-algebras
  \begin{itemize}
  \item is \emph{elementary} if it is among $x,\top,x_{1}\land
    x_{2},\bot,x_{1}\vee x_{2},\arrow[\sigma]X_{\sigma}$. 
  \end{itemize}
  With respect to two sets of variables $X$ and $Y$, we say that a
  term of the theory of modal $\mu$-algebras
  \begin{itemize}
  \item is \emph{simple} if it is a distributive combination of terms
    of the form $\bigwedge Y' \land
    \spcon[\emptyset,\Sigma]\set{D_{\sigma}}$, where $Y' \subseteq Y$
    and each $d \in D_{\sigma}$ is a distributive term on the
    variables in $X$,
  \item is \emph{disjunctive-simple} if it is a join of terms of the
    form $\bigwedge Y' \land
    \spcon[\emptyset,\Sigma]\set{D_{\sigma}}$, where $Y' \subseteq Y$
    and each $d \in D_{\sigma}$ is a join of a set of variables in
    $X$: $d = \bigvee X'$ with $X' \subseteq X$.
  \end{itemize} 
  \label{def:simpleF}
  For a $\mu$-algebra $\Alg{A}$ we say that a map $F = \langle
  F_{x}\rangle_{x \in X}: \Alg{A}^{X}\times \Alg{A}^{Y} \rTo
  \Alg{A}^{X}$ is \emph{elementary} (resp. \emph{simple}, resp.
  \emph{disjunctive-simple} w.r.t. $X$ and $Y$) if each component
  $F_{x} : \Alg{A}^{X}\times \Alg{A}^{Y} \rTo \Alg{A}$ is the
  interpretation of an elementary (resp. simple, resp.
  disjunctive-simple w.r.t. $X$ and $Y$) term.
\end{Definition}

\subsection*{Disjunctive-simple systems}
\begin{Proposition}
    Let $\Alg{F}$ be a free modal $\mu$-algebra, $G : \mathcal{F}^{X}
  \times \mathcal{F}^{Y}\rTo \mathcal{F}^{X}$ be a disjunctive-simple
  map, and let $k \in \mathcal{F}^{Y}$. Then $G_{k} : \mathcal{F}^{X}
  \rTo \mathcal{F}^{X}$ is a $\Opens$-adjoint of finite type.
\end{Proposition}
\begin{proof}  
  Proposition \ref{prop:list} and Lemma \ref{lemma:specconj} imply
  that for each $x \in X$ the $x$ component of $G_{k}$ -- which
  we shall denote $G_{x}$ abusing notation -- is a $\Opens$-adjoint.
  Item \ref{item:prods} in Proposition \ref{prop:list} then imply that
  $G_{k}$ is a $\Opens$-adjoint.  Thus we are mainly concerned with
  arguing that $G_{k}$ has finite type, and in view of Lemma
  \ref{lemma:uniformfinitary} and Lemma \ref{lemma:uniformfamily}, it
  will be enough to show that the $G_{x}$ form a uniform 
  family of finite type.
  
  Each $G_{x}$ has the form $\bigvee_{i \in I_{x}} k_{x,i} \land
  \spcon[\emptyset,\Sigma_{x,i}] \set{D_{\sigma}}$, where, for each $i
  \in I_{x}$, $k_{x,i}$ is a constant element of the free modal
  $\mu$-algebra $\mathcal{F}$, and for each $\sigma \in \Sigma_{x,i}$
  and $d \in D_{\sigma}$ $d = \bigvee X'$.  We use Lemma
  \ref{lemma:fincomp} and prove that the family
  \begin{align*}
    \mathcal{F} = \bigcup \set{\myjoin[Z]\,|\,Z \text{ is finite}}
    \cup \bigcup_{x \in X} \set{k_{x,i}\land z \,|\, i \in I_{x}} \cup
    \set{\spcon[\emptyset,\Sigma_{x,i}] \set{X_{\sigma}}}
  \end{align*}
  is uniform of finite type.  Here $\myjoin[Z] :
  \mathcal{F}^{Z} \rTo \mathcal{F}$ is the join operation of arity
  $Z$.
  
  For each constant $k_{x,i}$, choose a term $t_{x,i}$ representing
  the element $\neg k_{x,i}$.  Let now $Q_{0}$ be a finite subset
  of $\mathcal{F}^{X}$ and, for each $q \in Q_{0}$, let $s_{q}$ be a
  term representing $q$. 
  
  Let $FL$ be the Fisher-Ladner closure of the terms $t_{x,i}$ and
  $s_{q}$, it is well known \cite{kozen} that $FL$ is a finite set
  which, by its definition, comprises all the subterms of $t_{x,i}$
  and $s_{q}$.  
  
  Let $\overline{FL} \subseteq \mathcal{F}$ be the set of
  interpretations of terms in $FL$ in the $\mu$-algebra $\mathcal{F}$,
  and let $D$ be the distributive lattice generated by
  $\overline{FL}$. It will be useful to think of $D$ as the meet
  closure of the join closure of $\overline{FL}$.  We need to prove
  that $f \in \mathcal{F}$, $d \in D$, and $c \in \Cvrs{f}{d}$ with $c
  \in \mathcal{F}^{Z}$, imply $c(z) \in D$ for each $z \in Z$.
  
  Observe that if $c \in \Cvrs{\myjoin[Z]}{d}$, then $c$ is the vector
  with $d$ at each projection.

  If $f(z) = k_{x,i} \land z$, then $\Cvrs{f}{d} = \set{\neg k_{x,i}
    \vee d} \subseteq D$.
  
  Let $f = \spcon[\emptyset,\Sigma]$ and consider $d \in D$: since a
  cover in $\Cvrs{f}{d}$ is a meet of covers in $\Cvrs{f}{d_{i}}$
  where each $d_{i}$ belongs to the join closure of $\overline{FL}$,
  we can assume that $d$ is in the join closure of $\overline{FL}$,
  that is, $d$ is the interpretation of a term of the form $t_{1}\vee
  \ldots \vee t_{n}$ with $t_{i} \in FL$.
  Lemma \ref{lemma:specconj} shows that a $c \in
  \Cvrs{\spcon[\emptyset,\Sigma]}{d}$ is a meet of elements $c_{j}$,
  where each projection of a vector $c_{j}$ is either $\top$ or a join
  of subterms of $t_{1}\ldots t_{n}$, hence it  belongs to $D$.
\end{proof}

Lemma
\ref{lemma:finitaryconstructive} and the previous Lemma imply:
\begin{Corollary}
  \label{cor:simpleconstr} 
  The least prefixed point of a disjunctive-simple system $G :
  \mathcal{F}^{X}\times \mathcal{F}^{Y} \rTo \mathcal{F}^{X}$ is
  constructive:
  \begin{align*}
    \mu_{X}.G & = \bigvee_{n \geq 0} G_{v}^{n}(\bot)\,.
  \end{align*}
  for each $v \in \mathcal{F}^{Y}$.
\end{Corollary}

\subsection*{From disjunctive-simple to simple systems}
Our next goal is to transfer constructiveness from a disjunctive-simple
$G$ to a simple $F$. The main tool is the following Lemma:
\begin{Lemma}
  \label{lemma:transfer}
  Consider a commuting  diagram of posets with bottom
  $$
  \mydiagram[6em]{
    [](!s{L}{M}{L}{M}
    {1}{1},
    !a{^{i}}
    {^{f}}
    {^{g}}
    {^{i}}
    )
    "4":@/_1.5em/"3"_{\pi}
  }
  $$
  where $i$ is split by an order preserving $\pi$, $\pi \circ i =
  \id_{L}$.  Let $\alpha$ be a limit ordinal and suppose that (i) for
  $\beta < \alpha$, $f^{\beta}(\bot)$ and $g^{\beta}(\bot)$ exist and
  $i(f^{\beta}(\bot)) = g^{\beta}(\bot)$, (ii) the approximant
  $g^{\alpha}(\bot)$ exists. Then the approximant $f^{\alpha}(\bot)$
  exists as well and is equal to $\pi(g^{\alpha}(\bot))$.  If moreover
  $i$ is continuous, then $i(f^{\alpha}(\bot)) = g^{\alpha}(\bot)$.
\end{Lemma}
\begin{proof}
  Let $\alpha$ be an ordinal satisfying the hypothesis, we are going
  to argue that $\pi(g^{\alpha}(\bot)) = \bigvee_{\beta <\alpha}
  f^{\beta}(\bot)$. Let us begin supposing that, for some $l \in L$
  and every $\beta < \alpha$, $f^{\beta}(\bot) \leq l$.  Apply $i$
  to these relations and deduce that $g^{\beta}(\bot) \leq i(l)$ for
  $\beta < \alpha$, hence $g^{\alpha}(\bot) \leq i(l)$;  apply
  $\pi$ and deduce $\pi(g^{\alpha}(\bot)) \leq l$.  Conversely,
  apply $\pi$ to $i(f^{\beta}(\bot)) =
  g^{\beta}(\bot) \leq g^{\alpha}(\bot)$ to deduce $f^{\beta}(\bot) \leq \pi(g^{\alpha}(\bot))$ for $\beta < \alpha$.

  If moreover $i$ is continuous, then:
  $$
  i(f^{\alpha}(\bot)) 
  \;= \;
  i(\bigvee_{\beta < \alpha} f^{\beta}(\bot))
  \;=\;
  \bigvee_{\beta < \alpha} i(f^{\beta}(\bot))
  \;=\;
  \bigvee_{\beta < \alpha} g^{\beta}(\bot)
  \;= \;g^{\alpha}(\bot)\,.
   $$  
  
\end{proof}
We shall make use of the Lemma as follows. For a finite set of
variables $X$, let $\ppowerset(X)$ be the set of nonempty subsets of
$X$. For each $S \in \ppowerset(X)$, the map $i_{S}: L^{X} \rTo L$,
defined by
\begin{align*}
  i_{S}(x) & = \bigwedge_{j \in S} x_{j}\,,
\end{align*}
is continuous.  These maps, collected together, define  a
continuous map 
\begin{align*}
  i = \langle i_{S} \rangle_{S \in P_{+}(X)}& : L^{X} \rTo
  L^{\ppowerset(X)}\,
\end{align*}
which moreover preserves the bottom element.  For each $x \in X$
there is a projection onto the singleton set $\prj_{\set{x}} :
L^{\ppowerset(X)} \rTo L$.  These projections, collected into a common
projection $\prj = \langle \prj_{\set{x}}\rangle_{x \in X}:
L^{\ppowerset(X)} \rTo L^{X}$, split $i$: $\prj \circ i =
\id_{\mathcal{F}^{X}}$. Thus we shall prove:
\begin{Proposition}
  \label{prop:subsetconstruct}
  For each simple $F : \mathcal{F}^{X}\times \mathcal{F}^{Y} \rTo
  \mathcal{F}^{X}$ there is a disjunctive-simple $G :
  \mathcal{F}^{\ppowerset(X)}\times \mathcal{F}^{Y} \rTo
  \mathcal{F}^{\ppowerset(X)}$ such that the diagram
  \begin{align}
    \label{diag:subsetc}    
    & 
    \mydiagram[6em]{
      [](!s{\mathcal{F}^{X}\times \mathcal{F}^{Y}}
      {\mathcal{F}^{\ppowerset(X)}\times \mathcal{F}^{Y}}
      {\mathcal{F}^{X}}{\mathcal{F}^{\ppowerset(X)}}
      {1}{1.5},
      !a{^{i\times \id_{\mathcal{F}^{Y}}}}
      {^{F}}
      {^{G}}
      {^{i}}
      )
    }
  \end{align}
  commutes.
\end{Proposition}
Together with Corollary \ref{cor:simpleconstr} and Lemma
\ref{lemma:transfer}, the Proposition implies:
\begin{Corollary}
  Let $F : \mathcal{F}^{X}\times \mathcal{F}^{Y} \rTo
  \mathcal{F}^{X}$ be simple and $v \in \mathcal{F}^{Y}$. Then
  \begin{align*}
    \mu_{X}.F_{v} & = \bigvee_{n \geq 0} F_{v}^{n}(\bot)\,.
  \end{align*}
\end{Corollary}
Diagram \eqref{diag:subsetc} commutes if for each nonempty
subset $S\subseteq X$ we can find a disjunctive-simple $G_{S} :
\mathcal{F}^{\ppowerset(X)} \rTo \mathcal{F}$ such that
\begin{align}
  \label{eq:comm}
  \bigwedge_{j \in S} F_{j}(x) 
  & = 
  G_{S}(\bigwedge_{j \in S_{1}} x_{j},\ldots ,\bigwedge_{j \in S_{2^{n}-1}}
  x_{j})  \,,
\end{align}
where $S_{1},\ldots ,S_{2^{n}-1}$ is the list of nonempty subsets of
$X$. We shall sketch the proof of Proposition
\ref{prop:subsetconstruct}, skipping on the details since its
structure strictly follows \cite[\S 9.4]{AN}.  To be coherent with
this monograph, we use $ar(t_{1},\ldots ,t_{n})$ for the set of
variables in $X$ appearing in the terms $t_{i}$, while $t[t_{x}/x]$
denotes a standard substitution applied to the term $t$.
\begin{Lemma}
  For every pair of disjunctive-simple terms $t_{1},t_{2}$ there
  exists a disjunctive-simple term $t_{3}$ and $\phi :
  ar(t_{3})\rTo \ppowerset(ar(t_{1},t_{2}))$ such that the equation
  \begin{align*}
    t_{1} \land t_{2} & = t_{3}[\,\bigwedge \phi(x) /x\,]
  \end{align*}
  holds in every modal algebra.
\end{Lemma}
\begin{proof}
  Let $t_{1} = \bigvee_{i} \gamma_{1,i}$ and $t_{2} = \bigvee_{j}
  \gamma_{2,j}$, where the $\gamma_{k,l}$ have the form
  $\spcon[\Lambda,\Sigma]\set{D_{\sigma}}$ with $\Lambda \subseteq Y$
  and every $d \in D_{\sigma}$ a disjunction of variables from $X$.
  Clearly $t_{1} \land t_{2} = \bigvee_{i,j} \gamma_{1,i} \land
  \gamma_{2,j}$, thus it is enough to observe that
  \begin{align*}
    \spcon[\Lambda_{1},\Sigma_{1}]\set{D_{1,\sigma}} \land
    \spcon[\Lambda_{2},\Sigma_{2}]\set{D_{2,\sigma}} 
    & = (\,\spcon[\Lambda_{3},\Sigma_{3}]\set{D_{3,\sigma}}\,) 
    [\,\bigwedge \phi(x) /x\,]
  \end{align*} 
  for some some $\Lambda_{3},\Sigma_{3},D_{3,\sigma}$ and a $\phi$.
  To this goal, we observe that
  \begin{align*}
    \lefteqn[1cm]{\arrow[\sigma]D_{1,\sigma}
      \land \arrow[\sigma]D_{2,\sigma}}\\
    & = 
    \begin{cases}
      \arrow[\sigma]\emptyset, &   
      \text{if } D_{1,\sigma} = D_{2,\sigma} = \emptyset, \\
      \bot, &   
      \text{if } D_{i,\sigma}  = \emptyset \text{ for only one $i$},\\
      \lefteqn[4cm]{%
        \arrow[\sigma] 
        \set{d_{1} \land \bigvee D_{2,\sigma}
          ,  \bigvee D_{1,\sigma} \land d_{2}\,|\,d_{1} \in
          D_{1,\sigma},d_{2} \in D_{2,\sigma}},
      }\\
      &\text{otherwise}.
    \end{cases}
  \end{align*}
  We give the explicit definition of
  $\Lambda_{3},\Sigma_{3},D_{3,\sigma}$ and $\phi$, under the
  simplifying assumption that only the last case
  occurs.  We let $\Lambda_{3} = \Lambda_{1} \cup \Lambda_{2}$,
  $\Sigma_{3} = \Sigma_{1} \cup \Sigma_{2}$, $D_{3,\sigma} =
  D_{1,\sigma}$ for $\sigma \in \Sigma_{1}\setminus \Sigma_{2}$,
  $D_{3,\sigma} = D_{2,\sigma}$ for $\sigma \in \Sigma_{2}\setminus
  \Sigma_{3}$.  For $\sigma \in \Sigma_{1} \cap \Sigma_{2}$, we let
  \begin{align*}
    D_{3,\sigma} & =
    \set{
      \bigvee_{%
        \pile{x \in X \\ 
          \bigvee Z \in D_{2,\sigma} \\ 
          z \in Z}}
      w_{x,z}
      \;,\;
      \bigvee_{%
        \pile{y \in Y \\ 
          \bigvee Z \in D_{1,\sigma} \\ 
          z \in Z}}
      w_{z,y}
      \;|\;
      \bigvee X \in D_{1,\sigma}\,,\bigvee Y \in D_{2,\sigma}
      \,%
    }
  \end{align*}
  where $w_{x,z}$ and $w_{z,y}$ are new variables and let
  $\phi(w_{x,z}) = \set{x,z}$, $\phi(w_{z,y}) = \set{z,y}$.  
\end{proof}
The proof of Proposition \ref{prop:subsetconstruct} is then achieved
through the following steps:
\begin{enumerate}
\item It is shown that for every sequence of disjunctive-simple terms
  $t_{1},\ldots ,t_{n}$ there exists a disjunctive simple term $t_{0}$
  and $\phi : ar(t_{0})\rTo \ppowerset(ar(t_{1},\ldots ,t_{n}))$ such that
  \begin{align*}
    \bigwedge_{i=1,\ldots ,n} t_{i} & = t_{0}[\,\bigwedge 
    \phi(x)/x\,]
  \end{align*}
  is an equation of the theory of modal algebras.
\item It is shown that for each simple term $s$ there exists a
  disjunctive-simple term $d$ and a function $\phi: ar(d) \rTo
  \ppowerset(ar(s))$ such that
  \begin{align*}
    s & = d[\,\bigwedge \phi(x) /x\,]
  \end{align*}
  is an equation of the theory of modal algebras.
\end{enumerate}
Collecting together these properties, we obtain a rephrasing of
Proposition \ref{prop:subsetconstruct}:
\begin{enumerate}
\setcounter{enumi}{2}
\item For every sequence of  simple terms $s_{1},\ldots
  ,s_{n}$ there exists a disjunctive simple term $d$ and $\phi :
  ar(d)\rTo \ppowerset(ar(s_{1},\ldots ,s_{n}))$ such that
  \begin{align*}
    \bigwedge_{i=1,\ldots ,n} s_{i} & = d[\,\bigwedge 
      \phi(x)/x\,]
  \end{align*}
  is an equation of the theory of modal algebras.
\end{enumerate}

\subsection*{From simple to elementary systems}

Finally, we transfer constructiveness to elementary systems. To this
goal, we say that two systems $F: \mathcal{A}^{X} \times
\mathcal{A}^{Y} \rTo \mathcal{A}^{X}$ and $G: \mathcal{A}^{X} \times
\mathcal{A}^{Y} \rTo \mathcal{A}^{X}$ are \emph{equivalent} if for
each $v \in \mathcal{A}^{Y}$, the two chains of finite approximants
$\set{F^{n}_{v}(\bot)}_{n \geq 0}$ and $\set{G^{n}_{v}(\bot)}_{n \geq
  0}$ are cofinal into each other.  This means that for each $n \geq
0$ there exists $k \geq 0$ such that $F_{v}^{n}(\bot) \leq
G_{v}^{k}(\bot) $, and vice-versa.

\begin{Fact}
  If $F$ and $G$ are equivalent systems then $\bigvee_{n \geq
    0}F_{v}^{n}(\bot)$ exists if and only if $\bigvee_{n \geq 0} 
  G_{v}^{k}(\bot) $ exists, and in both cases they are equal.
\end{Fact}

We introduce now the notion of a \emph{guarded} system.  An occurrence
of a variable $x$ in a term $t$ is guarded if it is in the scope of a
modal operator. A term is guarded (w.r.t. $X$ and $Y$) if each
occurrence of a variable $x \in X$ in $t$ is guarded.  A system $F :
\mathcal{A}^{X} \times \mathcal{A}^{Y} \rTo \mathcal{A}^{X}$ is
guarded if each $F_{x}$ is the interpretation of a guarded term
(w.r.t. $X$ and $Y$).  The following Lemma is analogous to the well
known fact that every formula of the modal $\mu$-calculus is
equivalent to a guarded one \cite{kozen}. The reader may wish to
consult \cite[\S 9.2.4]{AN} as well.
\begin{Lemma}
  \label{lemma:epsilonelimination}
  For each elementary system $F$ there exists a guarded system $G$
  which is equivalent to $F$.
\end{Lemma}
\begin{proof}
  The following is a procedure -- analogous to $\epsilon$-transitions
  elimination in automata theory -- which eventually produces a system
  $G$, equivalent to a given system $F$, in which all the bound
  variables appear guarded.
  
  For a system $F$, let us define the graph of $\epsilon$-transitions:
  its nodes are the bound variables of $F$, and we say $x_{i}
  \rightarrow x_{j}$ if $x_{j}$ is not guarded in $F_{x_{i}}$.  Let us
  recall the  notion of distance between two nodes
  of this graph: if the two nodes $x_{i},x_{j}$ are connected by an
  $\epsilon$-path, then the distance between them is then minimum
  length of a path connecting them, and otherwise it is $\infty$.

  The procedure alternates among two kind of steps: elimination of
  loops, and reduction of cycles. 
  
  We can eliminate $\epsilon$-loops from $F$: disjunctive normal forms
  ensure that if $x_{i}$ is not guarded by a modal operator in
  $F_{x_{i}}$, then 
  \begin{align*}
    F_{x_{i}}(\ldots ,x_{i},\ldots ) & = (x_{i} \land f_{i}(\ldots
    ,x_{i},\ldots )) \vee g_{i}(\ldots ,x_{i},\ldots )
  \end{align*}
  for some terms $f_{i}$ and $g_{i}$ in which $x_{i}$ are guarded.
  We can modify $F$ to eliminate all the loops by the following
  rewrite of systems:
  $$
  \begin{array}{rcl}
    \begin{system}
      & \vdots \\
      x_{i} & = &  
      F_{x_{i}}(\ldots ,x_{i},\ldots ) \\
      & \vdots 
    \end{system}
    &
    \leadsto
    &
    \begin{system}
      & \vdots \\
      x_{i} & = &  
      g_{i}(\ldots ,x_{i},\ldots )\\
      & \vdots 
    \end{system}
  \end{array}
  $$
  Let $G$ be the system on the right, then it is easily checked
  that $F$ and $G$ have the same chains of approximants, hence they
  are equivalent.
  
  If $F$ does not contain loops, we can operate the following rewrite
  in which, for every pair of bound variables $x_{i}$ and $x_{j}$, if
  $x_{j}$ occurs unguarded in $F_{i}$ then it is substituted with
  $F_{j}$:
  $$
  \begin{array}{rcl}
    \begin{system}
      & \vdots \\
      x_{i} & = & F_{i}(\ldots ,x_{j},\ldots ) \\ 
      & \vdots
    \end{system}
     & \leadsto & 
    \begin{system}
      & \vdots \\
      x_{i} & = & F_{i}(\ldots ,F_{j}(\ldots ),\ldots ) \\ 
      & \vdots
    \end{system}
  \end{array}
  $$
  Observe that these rewrite reduce the distance between distinct
  connecetd nodes in the graph of $\epsilon$-transitions, thus the
  combined procedure terminates.  Let $G$ be the system on the right,
  then it is easily checked that $F_{v}^{n}(\bot) \leq G_{v}^{n}(\bot)
  \leq F_{v}^{2n}(\bot)$ for each $v \in {\cal A}^{Y}$, hence the two
  systems are equivalent.  
\end{proof}
The system $G$ obtained from the elementary $F$ by means of the
procedure described above need not to be a simple system.  On the
other hand, all the $G_{x}$ are terms of the theory of modal algebras
where all the variables in $X$ have modal depth at least $1$, so that
the system $G$ is quite similar to a simple one.  This means that by
adding new variables and cutting along substitutions we can
``unravel'' such a system $G$ to a simple system $H$. We only need to
justify these operations on systems.  To this goal, we modify the
previously proposed equivalence of systems.  Let $X \subseteq Z$, $F:
\Alg{A}^{X} \times \Alg{A}^{Y} \rTo \Alg{A}^{X}$, and $G:
\Alg{A}^{Z} \times \Alg{A}^{Y} \rTo \Alg{A}^{Z}$.  We say
that $G$ \emph{determines} $F$ iff the chains
$\set{F^{n}_{v}(\bot)}_{n \geq 0}$ and
$\set{\prj_{X}(G^{n}_{v}(\bot))}_{n \geq 0}$ are cofinal into each
other.
\begin{Fact}
  \label{fact:determining}
  Suppose that $G$ determines $F$, let $\Alg{A}$ be a modal
    $\mu$-algebra and  $v \in {\cal A}^{Y}$. If $\bigvee_{n \geq
    0}G^{n}_{v}(\bot)$ exists in $\Alg{A}^{Z}$, then $\bigvee_{n \geq
    0}F^{n}_{v}(\bot)$ exists in $\Alg{A}^{X}$ as well and is equal to
  $\prj_{X}(\bigvee_{n \geq 0}G^{n}_{v}(\bot))$.
\end{Fact}

\begin{Lemma}
  For each elementary system $F$ there exists a simple system
  $H$ which determines $F$.
\end{Lemma}
\begin{proof}
  We only sketch the proof. We apply the following kind of rewrite
  rules to the system $G$ obtained from $F$ by Lemma
  \ref{lemma:epsilonelimination}:
  \begin{align*}
    \begin{system}
      x_{1} & = & g(f(x_{1},x_{2},y),x_{1},x_{2},y) \\
      x_{2} & = & h(x_{1},x_{2},y) 
    \end{system}
    &    \leadsto
    \begin{system}
      x_{0} & = & f(x_{1},x_{2},y) \\
      x_{1} & = & g(f(x_{1},x_{2},y),x_{1},x_{2},y) \\
      x_{2} & = & h(x_{1},x_{2},y) 
    \end{system}\\[2mm]
    &    \leadsto
    \begin{system}
      x_{0} & = & f(x_{1},x_{2},y) \\
      x_{1} & = & g(x_{0},x_{1},x_{2},y) \\
      x_{2} & = & h(x_{1},x_{2},y) 
    \end{system}
\end{align*}
Let us call $G_{0}$, $G_{1}$, and $G_{2}$ the three systems in the
order.  Clearly $G_{1}$ determines $G_{0}$, while we have argued in
the proof of Lemma \ref{lemma:epsilonelimination} that $G_{2}$ is
equivalent to $G_{1}$. Hence $G_{2}$ determines $G_{0}$.  

Iteration of this rewriting produces a simple system $H$ determining
the original system $F$.  
\end{proof}

\begin{Proposition}
  Each elementary system is constructive on a free modal
  $\mu$-algebra.
\end{Proposition}
\begin{proof}
  In the previous subsection we have seen that simple systems are
  constructive on free modal $\mu$-algebras. In this subsection we
  have argued that given an elementary system there is a simple systems
  determining it, hence, on a free modal
  $\mu$-algebra, every elementary system is constructive by Fact
  \ref{fact:determining}.
\end{proof}


\section{$\Sigma_{1}$-operations are constructive}

 The valid equations
\begin{align*}
  \pos[\sigma]x & = \arrow[\sigma]\set{x,\top}
  &
  \nec[\sigma]x & =
  \arrow[\sigma]\set{x} \vee \arrow[\sigma]\emptyset
\end{align*}
show that all the operations of the theory of modal $\mu$-algebras are
definable from the Boolean algebra terms and the arrow terms
\eqref{eq:arrow}.  Accordingly we modify the definition
\eqref{eq:sigma1} of $\Sigma_{1}$-terms as follows:
$$
t \;=  \;x \,| \,\top \,| \,t \land t \,| \,
\bot \,|\, t \vee t \,|\, \arrow[\sigma] T
\,|\, \mu_{x}.t\,,
$$
where $x$ is a variable and $T$ is a set of previously defined
terms.  We remark that such a modification leaves invariant the class
of $\Sigma_{1}$-operations (i.e. interpretations of
$\Sigma_{1}$-terms).

The following concept is needed in the following:
\begin{Definition}
  We say that an order preserving map $f : L^{x}\times
  M^{y} \rTo L$ is \emph{regular} if it is continuous in each
  variable and constructive for the variable $x$. 
\end{Definition}
More generally, we shall say that an order preserving map $f : L^{X}
\rTo L$ is \emph{regular} if it is continuous and constructive in each
variable.  Recall that being constructive means that the approximant
$f_{v}^{\alpha}(\bot)$ exists for each $v \in M$ and each ordinal
$\alpha$.  It is easily seen that for a continuous $f$ existence of
the approximant $f_{v}^{\omega}(\bot) = \bigvee_{n \geq 0}
f^{n}_{v}(\bot)$ suffices for existence of all approximants. Hence,
when arguing that a continuous order preserving function is regular,
we shall only be concerned with existence of $f_{v}^{\omega}(\bot)$.
As an example, we have seen in the previous section that all the
elementary $G : \mathcal{F}^{X}\times \mathcal{F}^{Y} \rTo
\mathcal{F}^{X}$ are constructive on a free modal $\mu$-algebra
$\mathcal{F}$. Since each $G_{x} : \mathcal{F}^{X}\times
\mathcal{F}^{Y} \rTo \mathcal{F}$ is also continuous, $G$ is
continuous as well. Hence such elementary $G$ is regular on a free
modal $\mu$-algebra.  We want to transfer regularity, hence
constructiveness, to $\Sigma_{1}$-operations, for which we need to
consider them as solutions of elementary systems:
\begin{Lemma}
  \label{lemma:pi}
  For each $\Sigma_{1}$-operation $f: \mathcal{A}^{Y} \rTo
  \mathcal{A}$ there exists an elementary system $F: \mathcal{A}^{X}
  \times \mathcal{A}^{Y} \rTo \mathcal{A}^{X}$ and $x \in X$ such that
  $f = \prj_{x}\circ \mu_{X}.F$.
\end{Lemma}
\begin{proof} 
  The elementary systems are constructed by induction on the structure
  of $\Sigma_{1}$-terms.  For example, suppose $t = t_{1} \land t_{2}$
  and that $G_{i}:\Alg{A}^{X_{i}}\times \Alg{A}^{Y_{i}} \rTo
  \Alg{A}^{X_{i}}$ and $x_{i}$ have the property stated in the Lemma
  w.r.t. $t_{i}$, $i = 1,2$. We let $X = \set{x}\cup X_{1} \cup
  X_{2}$, $Y = Y_{1} \cup Y_{2}$, and $F$ is the system:
  $$
  \begin{system}
    x & = & x_{1} \land x_{2} \\
    x_{1} & = & G_{1,x_{1}}(\;\;\ldots\;\; ) \\
    & \vdots \\
    x_{2} & = & G_{2,x_{2}}(\;\;\ldots\;\; )\\
    & \vdots 
  \end{system}\,.
  $$
  Similar constructions work for $t = t_{1} \vee t_{2}$ and $t =
  \arrow[\sigma]T$. Suppose therefore that $t = \mu_{y}.t_{1}$. Assume
  that $G_{1}:\Alg{A}^{X_{1}}\times \Alg{A}^{Y_{1}} \rTo
  \Alg{A}^{X_{1}}$ and $x_{1}$ have the property stated in the Lemma
  w.r.t. $t_{1}$, then we let $X = X_{1} \cup \set{y}$, $Y = Y_{1}
  \setminus \set{y}$, and $F$ is the system:
  $$
  \begin{system}
    y & = & x_{1} \\
    x_{1} & = & G_{x_{1}}(\;\;\ldots\;\;) \\
    & \vdots 
  \end{system}\,.
  $$
  
\end{proof}

To achieve the proof of Claim \ref{claim:due} we need one more result,
which is a Beki\v{c}-like property for regular functions:
\begin{Proposition}
  \label{prop:bekicreg}
  Suppose that $F : L^{x} \times M^{y}\times N^{z} \rTo L$ is regular
  in $x$ and $G: L^{x} \times M^{y}\times N^{z} \rTo M$ is continuous.
  Then $\langle F,G\rangle : L^{x} \times M^{y}\times N^{z} \rTo L
  \times M$ is regular in $(x,y)$ if and only if $G \circ \langle
  \mu_{x}.F,\id_{M\times N}\rangle: M^{y} \times N^{z} \rTo M$ is
  regular in $y$.
\end{Proposition}
Using the Proposition above we can immediately state our goal:
\begin{Theorem}(cf.  Claim
\ref{claim:due}.)
  \label{theo:regular}
  Every $\Sigma_{1}$-operation $f : \mathcal{F}^{Z}\rTo \mathcal{F}$
  is regular on a free modal $\mu$-algebra $\mathcal{F}$, hence
  constructive.
\end{Theorem}
\begin{proof}
  By Lemma \ref{lemma:pi} $f = \prj_{x} \circ \mu_{X}.F$ for some
  elementary $F : \Alg{F}^{X} \times \Alg{F}^{Z} \rTo \Alg{F}^{X}$ and
  some $x \in X$.  Choose $y \in Z$ and observe that the system
  $$
  \langle F, \prj[X\cup Z]_{x}\rangle : \mathcal{F}^{X} \times
  \mathcal{F}^{y}\times \mathcal{F}^{Z \setminus \set{y}} \rTo
  \mathcal{F}^{X} \times \mathcal{F}^{y}$$
  is elementary hence
  regular. Since $\prj[X\cup Z]_{x}$ is continuous, we can use
  Proposition \ref{prop:bekicreg} with $G = \prj[X\cup Z]_{x}$ and
  deduce that $f = \prj[X]_{x} \circ \mu_{X}.F = \prj[X\cup Z]_{x} \circ
  \langle \mu_{X}.F,\id_{\mathcal{F}^{Z}} \rangle: \mathcal{F}^{Z} \rTo
  \mathcal{F}^{y}$ is regular for $y \in Z$.  
\end{proof}

Our main goal in the rest of the paper will be to prove Proposition
\ref{prop:bekicreg}. The next Lemma, needed often later, also
simplifies the statement of the Proposition.
\begin{Lemma}
  \label{lemma:continuous}
  If $g : L^{x} \times
  M^{y} \rTo M$ is regular, then $\mu_{y}.g : L^{x} \rTo M$
  is continuous.
\end{Lemma}
\begin{proof}
  \newcommand{\bigveeI}{\bigvee\!\! I}%
  Let $I$ be a directed set and suppose that $\bigvee I$ exists in
  $L$.  We argue first that $g^{m}_{\bigvee I}(\bot) = \bigvee_{i \in
    I} g^{m}_{i}(\bot)$ for all $m \geq 0$.  The relation trivially
  holds for $m = 0$.  Suppose it holds for $m$, then
  \begin{align*}
    g^{m+1}_{\bigveeI}(\bot) & = g(\bigvee I,g^{m}_{\bigveeI}(\bot)) 
    = g(\bigvee_{j \in I} j ,\bigvee_{i
      \in I}  g^{m}_{i}(\bot)) 
    \\
    & = \bigvee_{j \in I} \bigvee_{i
      \in I} g(j,g^{m}_{i}(\bot))
    = \bigvee_{i \in I} g(i,g^{m}_{i}(\bot)) = \bigvee_{i \in I} 
    g^{m+1}_{i}(\bot)\,.
    \intertext{Consequently, we obtain}
    (\mu_{y}.g)(\bigveeI)
    & =
    \bigvee_{m \geq 0}   g^{m}_{\bigveeI}(\bot) \\
    & =
    \bigvee_{m \geq 0}  \bigvee_{i \in I} g^{m}_{i}(\bot) 
    = \bigvee_{i \in I} \bigvee_{m \geq 0}  g^{m}_{i}(\bot)
    = \bigvee_{i \in I} (\mu_{y}.g)(i).
    \tag*{}
  \end{align*}  
\end{proof}%
The previous Lemma allows us to get rid of parameters in
Proposition \ref{prop:bekicreg} which we restate as follows:
\begin{Proposition}
  \label{prop:bekic}
  Consider $\langle f,g \rangle : L^{x} \times M^{y} \rTo L \times M$,
  where $f: L^{x} \times M^{y} \rTo L$ is continuous and $g: L^{x}
  \times M^{y} \rTo M$ is regular.  Then $\langle f,g \rangle$ is
  regular if and only if $f\circ \langle \id_{L},\mu_{y}.g\rangle:
  L^{x} \rTo L$ is regular.
\end{Proposition}
To prove the Proposition we shall fix a continuous $f$ and a regular
$g$.  We introduce an explicit notation for the approximants of
$\langle f,g\rangle$ and $h(x) = f(x,\mu_{y}.g(x,y))$:
\begin{align*}
  f_{0} & = \bot &
  g_{0} & = \bot \\
  f_{n +1} & = f(f_{n},g_{n})&
  g_{n+1} & = g(f_{n},g_{n})\,,
  \\
  \intertext{and}
  h_{0} & = \bot
  &
  i_{0} & = \bot \\
  h_{n + 1} & = f(h_{n},i_{n + 1})
  &
  i_{n+1} & = \mu_{y}.g(h_{n},y)\,.
\end{align*}
Using this notation we shall prove:
\begin{Proposition}
  \label{prop:bekic2}
  The sequenceces $\set{(f_{n},g_{n})}_{n \geq 0}$ and
  $\set{(h_{n},i_{n})}_{n \geq 0}$ have the same upper bounds.
\end{Proposition}

\noindent
\begin{proof}[Proof of Proposition \ref{prop:bekic}]
  Suppose $\langle f,g \rangle$ is regular. Since $g$ is regular,
  $\mu_{y}.g$ is continuous, by Lemma \ref{lemma:continuous}, hence
  $h(x) = f(x,\mu_{y}.g(x,y))$ is continuous, since $f$ is also
  continuous.  Thus we only need to show that $\bigvee_{n \geq 0}
  h^{n}(\bot) = \bigvee_{n \geq 0} h_{n}$ exists. Since $\langle f,g
  \rangle$ is regular $\bigvee_{n \geq 0} \langle f,g
  \rangle^{n}(\bot)$ exists, and by Proposition \ref{prop:bekic2}
  $\bigvee_{n \geq 0} \langle f,g \rangle^{n}(\bot) = \bigvee_{n \geq
    0} (f_{n},g_{n})= \bigvee_{n \geq 0} (h_{n},i_{n})$. Finally
  $\bigvee_{n\geq 0} h_{n}$ exists, since by continuity of projections
  it is equal to $\prj_{1}(\bigvee_{n \geq 0}
  (h_{n},i_{n}))$.
  
  Conversely, suppose that $h(x) = f(x,\mu_{y}.g(x,y))$ is regular.
  In particular we see that $\bigvee_{n\geq 0} h_{n}$ exists, and by
  Lemma \ref{lemma:continuous}, $\bigvee_{n\geq 0} i_{n} =
  \bigvee_{n\geq 0} \mu_{y}.g(h_{n},y) = \mu_{y}.g(\bigvee_{n\geq 0}
  h_{n},y)$ exists as well.  Clearly $\langle f,g\rangle$ is
  continuous and $\bigvee_{n \geq 0} \langle f,g \rangle^{n}(\bot) =
  \bigvee_{n \geq 0} (f_{n},g_{n}) = \bigvee_{n \geq 0} (h_{n},i_{n})$
  exists, by Proposition \ref{prop:bekic2}.  
\end{proof}

\noindent
\begin{proof}[Proof of Proposition \ref{prop:bekic2}]
  For a set $A$ we  let $M(A)$ be the set of upper bounds of $A$,
  i.e. $M(A) = \set{x\,|\,\forall a \in A\; a\leq x }$.
  
  Our first observation is that  $(f_{n},g_{n}) \leq
  (h_{n},i_{n})$, for all $n\geq 0$. 
  This is clearly true for $n = 0$, and if the 
  relation holds for $n$, then
  \begin{align*}
    f_{n+1} = f(f_{n},g_{n}) & \leq f(\nh_{n},\nii_{n})
    \leq f(\nh_{n},\nii_{n + 1}) = \nh_{n + 1}\,, \\
    g_{n+1} = g(f_{n},g_{n}) & \leq g(\nh_{n},\nii_{n}) \leq g(\nh_{n},\nii_{n+1})\\
    & = g(\nh_{n},\mu_{y}.g(\nh_{n},y)) 
    = \mu_{y}.g(\nh_{n},y) = \nii_{n +1} \,.
  \end{align*}
  Therefore we have $M(\set{(h_{n},i_{n})\,|\,n \geq 0}) \subseteq
  M(\set{(f_{n},g_{n})\,|\,n\geq 0})$.
  
  To prove the converse inclusion, we introduce a third sequence
  indexed by words of natural numbers:
  \begin{align*}
    \wi_{\epsilon} & = \bot
    &
    \wh_{\epsilon} & = \bot \\
    \wi_{wk} & = f(\wi_{w},g_{\wi_{w}}^{k}(\bot))
    &
    \wh_{wk} & = g_{\wi_{w}}^{k}(\bot)
    \,.
  \end{align*} 
  \begin{myclaim}
    The sequence $\set{(f_{n},g_{n})}_{n \geq 0}$ is cofinal into
    $\set{(l_{w},i_{w})}_{w \in \nnumbers^{\ast}}$:  for all $w \in
    \nnumbers^{\ast}$ there exists $n \in \nnumbers$ such that
    $(\wi_{w},\wh_{w}) \leq (f_{n},g_{n})$.
  \end{myclaim}%
  \begin{proof}[Proof of the Claim]
    Observe first that $g_{f_{n}}^{k}(\bot) \leq g_{n + k}$, for all
    $n,k \geq 0$.  This relation is trivial if $k = 0$, and supposing
    it holds for $k$, then
    \begin{align*}
      g^{k + 1}_{f_{n}}(\bot) 
      & = g(f_{n},g^{k}_{f_{n}}(\bot)) 
      \leq g(f_{n + k},g_{n + k}) = g_{n + k + 1}\,.
      \intertext{Clearly $(\wi_{\epsilon},\wh_{\epsilon}) \leq
        (f_{0},g_{0})$, and if $(\wi_{w},\wh_{w}) \leq
        (f_{n},g_{n})$, then} 
      \wh_{wk} & = g_{\wi_{w}}^{k}(\bot)
      \leq g_{f_{n}}^{k}(\bot) 
      \leq g_{n + k} \leq  g_{n + k + 1} \\
      \wi_{wk} & = f(\wi_{w},\wh_{wk})
      \leq f(f_{n},g_{n + k}) 
      \leq f(f_{n +k},g_{n + k}) = f_{n + k + 1}\,.
    \tag*{}
  \end{align*}
\end{proof}%
Consequenlty we have $M(\set{(f_{n},g_{n})\,|\,n\geq 0})\subseteq
M(\set{(l_{w},g_{w})\,|\,w \in \nnumbers^{\ast}})$.
  \begin{myclaim}
    The two relations
    \begin{align*}
      \nh_{n} & = 
      \bigvee_{w \in \nnumbers^{n}} 
      \wi_{w}
      & 
      \nii_{n} & = \bigvee_{w \in \nnumbers^{n}}
      \wh_{w} 
    \end{align*}
    hold.
  \end{myclaim}%
  \begin{proof}[Proof of the Claim]
    Observe first that if $w ,u \in \nnumbers^{n}$ and $w \leq u$,
    then $l_{w} \leq l_{u}$ and $m_{w} \leq m_{u}$. This is easily
    verified by induction on $n$. It follows that for $n$ fixed, the
    sets $\set{l_{w}}_{w \in \nnumbers^{n}}$ and $\set{m_{w}}_{w \in
      \nnumbers^{n}}$ are directed.
    
    The relations stated in the Claim trivially hold for $n = 0$.
    Suppose they hold for $n$. Then
    \begin{align*}
      \nii_{n + 1} & = \mu_{y}.g(\nh_{n},y) 
      = \mu_{y}.g(\bigvee_{w \in \nnumbers^{n}}\wi_{w},y) 
      = \bigvee_{w \in \nnumbers^{n}}
      \mu_{y}.g(\wi_{w},y)
    \\&
    = \bigvee_{w \in \nnumbers^{n}} \bigvee_{k\geq 0}
    g^{k}_{\wi_{w}}(\bot) 
    = \bigvee_{w \in \nnumbers^{n}} \bigvee_{k\geq 0}
    \wh_{wk} 
    = \bigvee_{u \in \nnumbers^{n +1}} 
    \wh_{u} 
    \intertext{and}
    \nh_{n + 1} & = f(\nh_{n},\mu_{y}.g(\nh_{n},y)) 
    = \bigvee_{w \in \nnumbers^{n}}
    f(\wi_{w},\mu_{y}.g(\wi_{w},y))    \\
    & = \bigvee_{w \in \nnumbers^{n}}
    f(\wi_{w},\bigvee_{k\geq 0}g^{k}_{\wi_{w}}(\bot))  
    = \bigvee_{w \in \nnumbers^{n}}\bigvee_{k\geq 0}
    f(\wi_{w},g^{k}_{\wi_{w}}(\bot))    
    = \bigvee_{u \in \nnumbers^{n+1}}
    \wi_{u}\,.
    \tag*{}
  \end{align*}
  \end{proof}%
  Consequenlty, $M(\set{(h_{n},i_{n})\,|\,n\geq 0})=
  M(\set{(l_{w},m_{w})\,|\,w \in \nnumbers^{\ast}})$ and, by the
  previous results, $M(\set{(h_{n},i_{n})\,|\,n\geq 0})=
  M(\set{(f_{n},g_{n})\,|\,n\geq 0})$.  This terminates the proof of
  Proposition \ref{prop:bekic2}.  
\end{proof}


\bibliographystyle{plain}
\bibliography{../biblio}

\end{document}